\documentclass[11pt,a4paper]{amsart}
\usepackage{mathrsfs}
\usepackage{amsfonts}
\usepackage[notref,notcite]{showkeys}
\usepackage{enumerate}

\usepackage{amsmath,amssymb,xspace,amsthm}
\newtheorem{theorem}{Theorem}[section]
\newtheorem{lemma}[theorem]{Lemma}%[section]
\newtheorem{definition}{Definition}%[section]
\numberwithin{equation}{section}

\def\Der{\operatorname{Der}}

\newcommand{\C}{\ensuremath{\mathbb C}\xspace}

\newcommand{\h}{\ensuremath{\mathfrak{h}}}

\newcommand{\Z}{\ensuremath{\mathbb{Z}}\xspace}

\newcommand{\W}{\ensuremath{\mathcal{W}}\xspace}

\renewcommand{\phi}{\varphi}

\renewcommand{\leq}{\leqslant}
\renewcommand{\geq}{\geqslant}

\def\sl{\mathfrak{sl}}

\def\A{\mathcal{A}}

\begin{document}
\title[Non-weight representations of Cartan type S Lie algebras] {Non-weight representations of Cartan type S Lie algebras}

%\date{}
\maketitle
\centerline{Juanjuan Zhang$^{1,2}$}
\begin{center}
{\small
$^1$College of Mathematics and Information Science, Hebei Normal (Teachers) University, Shijiazhuang, Hebei, 050016, PR China\\
$^2$School of Mathematics and Statistics, Qingdao University, Qingdao, Shandong, 266071, PR China\\}
\end{center}
\begin{abstract}
For the two Cartan type S subalgebras of the Witt algebra $\W_n$, called Lie algebras of divergence-zero vector fields, we determine all module structures on the universal enveloping algebra of their Cartan subalgebra $\h_n$. We also give all submodules of these modules.
\end{abstract}


\vskip 10pt
2000 Math. Subj. Class.: 17B10, 17B20, 17B65, 17B66, 17B68 

\vskip 10pt
\section{introduction}

\vskip 5pt
Throughout this paper, we denote by $\mathbb{Z}$, $\mathbb{Z}^{\ast}$, $\mathbb{Z}_+$, $\mathbb{N}$,
$\mathbb{C}$ and $\mathbb{C}^{\ast}$ the sets of  all integers, nonzero integers, non-negative integers,
positive integers, complex numbers and nonzero complex numbers, respectively. All vector
spaces and algebras are over $\C$. We denote by  $U(\mathfrak{a})$ the universal enveloping algebra of a Lie algebra $\mathfrak{a}$,
and by $\mathfrak{a}'$ the derived algebra.

Classification of simple modules is a natural question which arises when
studying the representation theory of Lie algebras. Simple modules are, in some
sense, “building blocks” for all other modules, and hence understanding them
is important. But, the problem to classify all simple modules over nontrivial Lie algebras has been considered impossible, see Diximier's remark in [B].
However, because of its importance, the problem of construction of new
families of modules attracted a lot of attention over the years. The most studied case
seems to be the Virasoro  algebra, where many different multi-parameter
families of simple modules were constructed by various authors; see, for example, [OW],
[GLZ], [LZ], [LLZ], [MZ1], [MZ2], [MW] and references therein.

Recently, a class of modules over $\mathfrak{sl}_{n+1}$ and Witt algebras were constructed in [N2] and [TZ1] independently. These modules have free actions of the Cartan subalgebra are studied intensively. More precisely, in [N2], Jonathan Nilsson classified the $U(\mathfrak{sl}_{n+1})$-modules whose restrictions to $U(\h)$ are free modules of rank 1, and determined their submodules ($\h$ is the standard Cartan subalgebra of $\mathfrak{sl}_{n+1}$ consisting of all traceless diagonal matrices).
In [N3],  Jonathan Nilsson investigated the category of $U(\h)$-free $\mathfrak{g}$-modules (where $\mathfrak{g}$ is a finite-dimensional simple complex Lie algebra with a fixed Cartan subalgebra $\h$), showing that $U(\h)$-free modules only could exist when $\mathfrak{g}$ is of type A or C, and classifying isomorphism classes of $U(\h)$-free modules of rank 1 in type C, which includes an explicit construction of new simple $\mathfrak{sp}(2n)$-modules.
In [TZ1], a class of $\W_n$-modules $\Omega(\Lambda_n,a)$ was defined, which turned to be (see [TZ2])
all the possible $\W_{n}$-modules which are free $U(\h_{n})$-modules of rank 1 ($\h_{n}$ is the fixed Cartan subalgebra of $\W_n$), where the simplicity of such modules was also completely decided.

Motivated by [N2] and [TZ2],  we intend to study such modules for the two Cartan type S Lie algebras $\tilde S_n$ and $\bar S_n$, also called Lie algebras of divergence zero vector fields, and to determine their submodule structures, too.

First, let's recall the algebras.
For any integer $n : 2\le n \le \infty$, let
$\A_n=\C[t_1^{\pm1},t_2^{\pm1},\cdots, t_n^{\pm1}]$  be the algebra of Laurent polynomials over $\C$ in $n$ variables $t_1,t_2,\cdots, t_n$.
Then $\W_n=\Der(\A_n)$ is the Witt algebra of rank $n$, elements of which are the vector fields:

$$\begin{aligned}
\W_n&=\bigoplus_{i=1}^n\A_n \partial_i\\
&=\left\{\sum_{i=1}^{n}f_{i}\partial_i \bigg| f_{i}\in\A_n,i=1,\cdot\cdot\cdot,n \right\}\\
&=\left\{\sum_{i=1}^{n}f_{i}\frac{\partial}{\partial t_i} \bigg| f_{i}\in\A_n,i=1,\cdot\cdot\cdot,n\right\}.
\end{aligned}$$
There is a standard
{\emph {Cartan subalgebra}} $\mathfrak{h}_{n}=\oplus_{i=1}^n\C \partial_i$ of $\W_n$,
where $\partial_i= t_i\frac {\partial}{\partial t_i},i=1,\cdot\cdot\cdot,n$. For $r=(r_{1},\cdot\cdot\cdot,r_{n})\in \Z^n$,
set $t^{r}=t_{1}^{r_{1}}t_{2}^{r_{2}}\cdot\cdot\cdot t_{n}^{r_{n}}$. For later use, we set $\mathbf{-1}=(-1, \cdot\cdot\cdot, -1)\in \Z^n$ and
$\mathbf{0}=(0, \cdot\cdot\cdot, 0)\in \Z^n$. The Lie brackets in $\W_n$ are:
$$[t^r\partial_i,t^s\partial_j]=s_it^{r+s}\partial_j-r_jt^{r+s}\partial_i, \,\,\forall \,\, i,j=1,2,\cdots, n; r,s\in \Z^n.$$

Let $\text{div}$ be the divergence, defined by ([DZ] Page 147)
$$\text{div}\left(\sum_{i=1}^{n}f_{i}\partial_i\right)=\sum_{i=1}^{n}\partial_i(f_{i}),$$
$\text{div}$ is a derivation of $\W_n$ with values in $\A_n$:
$$\text{div}[u,v]=u\cdot\text{div}(v)+v\cdot\text{div}(u), u, v\in\W_n$$
where the dot indicates the natural action of $\W_n$ on $\A_n$. In view of the above equation,
$$\begin {aligned}
\tilde{S}_{n}=\ker(\text{div})&=\left\{\sum_{i=1}^{n}f_{i}\partial_i\bigg|f_{i}\in\A_n,i=1,\cdot\cdot\cdot,n,\sum_{i=1}^n\partial_i(f_{i})=0\right\}\\
&=\text{Span}_{\mathbb{C}}\{\partial_i,r_{j}t^{r}\partial_i-r_{i}t^{r}\partial_j\mid i,j\in\{1,\cdot\cdot\cdot,n\},r\in \Z^n\},
\end{aligned}$$
is a subalgebra of $\W_n$ ([DZ] Page 147, where $\tilde{S}=\tilde{S}_{n}$). This algebra is called the Lie algebra of divergence zero vector fields.

The derived subalgebra of $\tilde{S}_{n}$:
$$\tilde{S}_{n}'=\text{Span}_{\mathbb{C}}\{r_{j}t^{r}\partial_i-r_{i}t^{r}\partial_j\mid i,j\in\{1,\cdot\cdot\cdot,n\},r\in \Z^n\}.$$

The algebra $\tilde{S}_{n}$ is not simple, but its derived subalgebra $\tilde{S}_{n}'$ is simple, assuming only that $n \geq 2$ ([DZ] Page 145).

We have the classical divergence defined as follows:
$$\text{Div}\left(\sum_{i=1}^{n}f_{i}\frac{\partial}{\partial t_i}\right)=\sum_{i=1}^{n}\frac{\partial f_{i}}{\partial t_i}\ ([DZ]\ Page\ 147).$$
Then we have the subalgebras
$$\tilde{S}_{cl}=\ker(\text{Div})=t^{\mathbf{-1}}\tilde{S}_{n}\ ([DZ]\ Page\ 149)$$
and ([DZ] Proposition 3.2)
$$\begin{aligned}
\tilde{S}_{cl}'&=(t^{\mathbf{-1}}\widetilde{S}_{n})'=t^{\mathbf{-1}}\tilde{S}_{n}'\\
&=\text{Span}_{\mathbb{C}}\{(r_{j}+1)t^{r}\partial_i-(r_{i}+1)t^{r}\partial_j\mid i,j\in\{1,\cdot\cdot\cdot,n\},r\in \Z^n\}.
\end{aligned}$$
We know that $\tilde{S}_{cl}'=t^{\mathbf{-1}}\tilde{S}_{n}'$ is simple when $n\geq3$ ([DZ] Theorem 3.3); when $n=2$, the Lie algebras $t^{\mathbf{-1}}\tilde{S}_{2}'$ is not simple,
but its derived algebra $(t^{\mathbf{-1}}\tilde{S}_{2}')'$ is simple and has codimension 1 in $t^{\mathbf{-1}}\tilde{S}_{2}'$.
We refer to the simple Lie algebras $(t^{\mathbf{-1}}\tilde{S}_{n}')'$  as the Lie algebras of generalized Cartan type S
([DZ]  Page 145, Page 150-152). Set
$$\bar{S}_{n}=\mathfrak{h}_{n}\ltimes (t^{\mathbf{-1}}\tilde{S}_{n}')'.$$
Specifically, for $n\geq3$,
$$\bar{S}_{n}=\text{Span}_{\mathbb{C}}\{\partial_i, (r_{j}+1)t^{r}\partial_i-(r_{i}+1)t^{r}\partial_j \mid i,j\in\{1,\cdot\cdot\cdot,n\}, r\in \mathbb{Z}^{n} \},$$
and for $n=2$,
$$\bar{S}_{2}=\text{Span}_{\mathbb{C}}\{\partial_1, \partial_2, (r_{2}+1)t^{r}\partial_1-(r_{1}+1)t^{r}\partial_2 \mid r\in \mathbb{Z}^{2},r\neq(-2,-2)\}.$$
Then $\bar{S}_{n}'=t^{\mathbf{-1}}\tilde{S}_{n}'$ for $n\geq3$ and
$$\bar{S}_{2}'=(t^{\mathbf{-1}}\tilde{S}_{2}')'=\text{Span}_{\mathbb{C}}\{ (r_{2}+1)t^{r}\partial_1-(r_{1}+1)t^{r}\partial_2 \mid r\in \mathbb{Z}^{2},r\neq(-2,-2)\}$$
are simple.

Note that $\tilde{S}_{n}$ and $\bar{S}_{n}$ are subalgebras of $\W_n$, and they share the {\emph {Cartan subalgebra}} $\mathfrak{h}_{n}=\oplus_{i=1}^n\C \partial_i$.\\

Next, let's introduce the modules defined in [TZ1]. Let $2\le n\le\infty, a\in \C, \Lambda_n=(\lambda_1,\lambda_2,\cdots,\lambda_n)\in (\C^*)^n$.
Denote by $\Omega(\Lambda_n,a)=\C[x_1,x_2,\cdots, x_n]$ the polynomial  algebra over $\C$
in commuting indeterminates $x_1,x_2,\cdots, x_n.$
The action of $\W_n$ on $\Omega(\Lambda_n,a)$ is defined by
\begin{equation}
t^k\partial_i \cdot f(x_1,\cdots,x_n)=\Lambda_n^k(x_i-k_i(a+1))f(x_1-k_1,\cdots,x_n-k_n),
\end{equation}
where $k=(k_1,k_2,\cdots,k_n)\in \Z^n,f(x_1,\cdots,x_n)\in \C[x_1,x_2,\cdots, x_n]$,
$\Lambda_n^k=\lambda_1^{k_1}\lambda_2^{k_2}\cdots\lambda_n^{k_n}$,
$i=1,2,\cdots,n$.

Since $\tilde{S}_{n}$ and $\bar{S}_{n}$ are subalgebras of $\W_n$, $\Omega(\Lambda_n,a)$ can be seen as $\tilde{S}_{n}$- and $\bar{S}_{n}$-modules, which will still be denoted by $\Omega(\Lambda_n,a)$ without causing confusion.

For any $n$, $2\le n\le\infty$, define a class of $\tilde{S}_{n}$-modules as follows:
\begin{definition} Let $\Lambda_n=(\lambda_1,\cdots,\lambda_n)\in (\C^*)^n, \alpha=(\alpha_{1},\cdot\cdot\cdot,\alpha_{n})\in \C^{n}$,
define the action of $\tilde{S}_{n}$ on $\C[\partial_1,\cdot\cdot\cdot,\partial_n]$ as follows:
\begin{equation}
\begin{aligned}
&(r_{j}t^{r}\partial_i-r_{i}t^{r}\partial_j)f(\partial_1,\cdot\cdot\cdot,\partial_n)\\
&=\Lambda_n^{r}(r_{j}(\partial_i+\alpha_{i})-r_{i}(\partial_j+\alpha_{j}))
f(\partial_1-r_{1},\cdot\cdot\cdot,\partial_n-r_{n}),
\end{aligned}
\end{equation}
where $r=(r_1,\cdots,r_n)\in\Z^n$,
$f(\partial_1,\cdot\cdot\cdot\partial_n)\in\C[\partial_1,\cdot\cdot\cdot,\partial_n]$, and $\partial_i (i=1,2,\cdots,n)$ act by left multiplication.
We denote this $\tilde{S}_{n}$-module by  $\Omega(\Lambda_n,\alpha)$.
\end{definition}

It is straightforward to verify that the above action makes  $\C[\partial_1,\cdot\cdot\cdot,\partial_n]$ into a $\tilde{S}_{n}$-module.
Very surprisingly, only the $\tilde S_n$-modules  $\Omega(\Lambda_n,\mathbf{0})$  is the restriction of  $\Omega(\Lambda_n,a)$ from $\W_n$,
and these $\tilde S_n$-modules are parameterized by additional n constants.
When $n=2$, the $\tilde{S}_{2}$-modules, $\Omega(\Lambda_2,\alpha)$ are also denoted by $\Omega(\lambda_{1},\lambda_{2},\alpha_{1},\alpha_{2})$,
with $\Lambda_2=(\lambda_1, \lambda_2)\in (\C^*)^2$ and $\alpha=(\alpha_{1}, \alpha_{2})\in \C^{2}$.

For any $n$, $2\le n\le\infty$, define a class of $\bar{S}_{n}$-modules as follows:

\begin{definition} Let $\Lambda_n=(\lambda_1,\cdots,\lambda_n)\in (\C^*)^n, \kappa\in \C$,
define the action of $\bar{S}_{n}$ (respectively $\bar{S}_{2}$) on $\C[\partial_1,\cdot\cdot\cdot,\partial_n]$ as follows:
\begin{equation}
\begin{aligned}
&((r_{j}+1)t^{r}\partial_i-(r_{i}+1)t^{r}\partial_j)f(\partial_1,\cdot\cdot\cdot,\partial_n)\\
&=\Lambda_n^{r}((r_{j}+1)\partial_i-(r_{i}+1)\partial_j+(r_{i}-r_{j})\kappa)f(\partial_1-r_{1},\cdot\cdot\cdot,\partial_n-r_{n}),
\end{aligned}
\end{equation}
where $r=(r_1,\cdots,r_n)\in\Z^n$,
$f(\partial_1,\cdot\cdot\cdot\partial_n)\in\C[\partial_1,\cdot\cdot\cdot,\partial_n]$, and $\partial_i (i=1,2,\cdots,n)$ act by left multiplication.
We denote this $\bar{S}_{n}$ (respectively $\bar{S}_{2}$)-module by $\Omega(\Lambda_n,\kappa)$.
\end{definition}
The $\bar{S}_{n}$-modules $\Omega(\Lambda_n,\kappa)$ is only the restriction of  $\Omega(\Lambda_n,a)$ from $\W_n$ with $\kappa=-(a+1)$.
Unlike the above case, these modules do not have more parameters.

The present paper is organized as follows.
In Section 2, we first work out the $\tilde{S}_{2}$-modules that are free $U(\h_2)$-modules of rank 1 (Theorem 2.6); then we picture clearly the $\tilde{S}_{n}$-modules that are free $U(\h_n)$-modules of rank 1 (Theorem 2.13). These modules are exactly $\Omega(\Lambda_n,\alpha)$ for some $\Lambda_n=(\lambda_1,\cdots,\lambda_n)\in (\C^*)^n$, $\alpha=(\alpha_{1},\cdot\cdot\cdot,\alpha_{n})\in \C^{n}$.

In Section 3, we first work out the $\bar{S}_{2}$-module structures  on $U(\h_2)$ (Theorem 3.6); then we find out the $\bar{S}_{n}$-modules that are free $U(\h_n)$-modules of rank 1 (Theorem 3.12). These modules are exactly $\Omega(\Lambda_n,\kappa)$ for some $\Lambda_n=(\lambda_1,\cdots,\lambda_n)\in (\C^*)^n$, $\kappa\in \C$.

\section{$\tilde{S}_{n}$-module structures  on $U(\h_n)$}

\vskip 5pt

\subsection{The Virasoro-like algebra}

The Virasoro-like algebra $\mathscr{L}$ is a Lie algebra with basis ([WT] Page 4163)

$$\{d_{i,j},h_{1},h_{2} \mid (i,j)\in\mathbb{Z}^{2}\setminus\{(0,0)\}\}$$
and brackets
$$[d_{i,j},d_{k,l}]=(jk-il)d_{i+k,j+l};$$
$$[h_{1},d_{i,j}]=id_{i,j},[h_{2},d_{i,j}]=jd_{i,j},[h_{1},h_{2}]=0,$$
and $d_{0,0}=0$ understandably. The derived subalgebra $\mathscr{L}'$ has the basis $\{d_{i,j}\mid(i,j)\in\mathbb{Z}^{2}\setminus\{(0,0)\}\}$, which is simple ([K] Page 3762, [Z] Page 506, where $V=\mathscr{L}'$) and can be generated by $d_{1,0},d_{-1,0},d_{0,1},d_{0,-1}$([K] Page 3761). The Cartan subalgebra is $\mathfrak{h}_2=\mathbb{C}\mathbb{c}h_{1}\oplus\mathbb{C}\mathbb{c}h_{2}$, and $\mathscr{L}=\mathscr{L}'\oplus \mathfrak{h}_2$.

In $U(\h_2)=\mathbb{C}[h_{1},h_{2}]$, for $i=1,2$, define ([N2] Page 298, 299)
$$\text {deg}_ih_{1}^{d_{1}}h_{2}^{d_{2}}=d_{i},\text {deg}_i0=-1$$
and commutative algebra automorphisms
$$\sigma_i:  U(\h_2)\rightarrow U(\h_2)$$ $$ h_k\mapsto h_k-\delta_{ik}$$

For $f\in\mathbb{C}[h_{1},h_{2}]$,
$$\sigma_1(f(h_{1},h_{2}))=f(h_1-1,h_{2}),\sigma_2(f(h_{1},h_{2}))=f(h_{1},h_{2}-1)$$
and
$$\text {deg}_{i}(\sigma_{i}(f)-f)=(\text {deg}_{i}f)-1.$$

Note that definitions of $\text {deg}_{i}$ and $\sigma_{i}$ can be extended uniquely to $\mathbb{C}[h_{1},h_{2},\cdot\cdot\cdot, h_{n}]$.

The Virasoro-like algebra $\mathscr{L}$ is a subalgebra of the Witt algebra $\W_2$, and we have that:

$$\begin{aligned}
\phi:&\mathscr{L}\rightarrow\W_2\\
&h_{1}\mapsto\partial_1\\
&h_{2}\mapsto\partial_2\\
&d_{1,0}\mapsto t_{1}\partial_2\\
&d_{-1,0}\mapsto -t_{1}^{-1}\partial_2\\
&d_{0, 1}\mapsto t_{2}\partial_1\\
&d_{0, -1}\mapsto-t_{2}^{-1}\partial_1\\
&(-1)^{i}d_{i,j}\mapsto jt_{1}^{i}t_{2}^{j}\partial_1-it_{1}^{i}t_{2}^{j}\partial_2,
\end{aligned}$$
where $i,j\in\mathbb{Z}, (i,j)\neq(0,0)$, $\partial_k=t_k\frac{\partial}{\partial t_k}, k=1,2$, and $\phi$ is an injective homomorphism of Lie algebras. Clearly, $\phi(\mathscr{L})=\tilde{S}_{2}$, $\phi(\mathscr{L}')=\tilde{S}_{2}'$, i.e., $\mathscr{L}\cong\tilde{S}_{2}$, $\mathscr{L}'\cong\tilde{S}_{2}'$, and we'll not distinguish  them in the following.\\

\subsection{$\mathscr{L}$-module structures on $U(\h_2)$}

Let $M=U(\h_2)\cdot1\cong\mathbb{C}[h_{1},h_{2}]$ be a $\mathscr{L}$-module that is a free $U(\h_2)$-module of rank 1, on which $\mathscr{L}'$ acts nontrivially.
\begin{lemma} For any $f(h_{1},h_{2})\cdot1\in M,(i,j)\in\mathbb{Z}^{2}\setminus\{(0,0)\}$,
\begin{equation}d_{i,j}\cdot f(h_{1},h_{2})1
=f(h_{1}-i,h_{2}-j)(d_{i,j}\cdot1),
\end{equation} and $d_{i,j}\cdot1\ne0$.
\end{lemma}
\begin{proof} Formula (2.1) follows from
$$d_{i,j}h_{1}=(h_{1}-i)d_{i,j},d_{i,j}h_{2}=(h_{2}-j)d_{i,j}\in U(\mathscr{L}).$$

Now assume that $d_{i,j}\cdot1=0$ for some $(i,j)\ne (0,0)$,
then  $d_{i,j}\cdot M=0$, i.e., ann$(M)\ne0$. Since  $\mathscr{L}'$ is  simple
 and ann$(M)$ is a nonzero ideal of $\mathscr{L}'$,
this forces ann$(M)=\mathscr{L}'$, a contradiction.
\end{proof}\

Since the Lie algebra $\mathscr{L}'$ is generated by $d_{1,0},d_{-1,0},d_{0,1},d_{0,-1}$, by Lemma 2.1, the module structures on $M$ are determined by $$d_{1,0}\cdot1, d_{-1,0}\cdot1, d_{0,1}\cdot1, d_{0,-1}\cdot1.$$
There are nonzero polynomials
$$f_{\pm1,0}(h_{1},h_{2}), f_{0,\pm1}(h_{1},h_{2})\in\mathbb{C}[h_{1},h_{2}]$$
such that
$$d_{\pm1,0}\cdot1=f_{\pm1,0}(h_{1},h_{2}), d_{0,\pm1}\cdot1=f_{0,\pm1}(h_{1},h_{2}).$$

\begin{lemma}
$$f_{\pm1,0}(h_{1},h_{2})=f_{\pm1,0}(h_{2})\in\mathbb{C}[h_{2}],f_{0,\pm1}(h_{1},h_{2})=f_{0,\pm1}(h_{1})\in\mathbb{C}[h_{1}].$$
\end{lemma}

\begin{proof}
Since $[d_{1,0},d_{-1,0}]\cdot1=0$, i.e.,
$$\begin{aligned}
&d_{1,0}\cdot f_{-1,0}(h_{1},h_{2})1=d_{-1,0}\cdot f_{1,0}(h_{1},h_{2})1,\\
&f_{-1,0}(h_{1}-1,h_{2})f_{1,0}(h_{1},h_{2})=f_{1,0}(h_{1}+1,h_{2})f_{-1,0}(h_{1},h_{2}),\\
&\sigma_{1}(f_{1,0}(h_{1}+1,h_{2})f_{-1,0}(h_{1},h_{2}))=f_{1,0}(h_{1}+1,h_{2})f_{-1,0}(h_{1},h_{2}),
\end{aligned}$$
so $f_{1,0}(h_{1}+1,h_{2})f_{-1,0}(h_{1},h_{2})\in\mathbb{C}[h_{2}]$, and $f_{\pm1,0}(h_{1},h_{2})\in\mathbb{C}[h_{2}]$.
Write simply $ f_{\pm1,0}(h_{1},h_{2})=f_{\pm1,0}(h_{2})$. By similar arguments, $f_{0,\pm1}(h_{1},h_{2})=f_{0,\pm1}(h_{1})\in\mathbb{C}[h_{1}]$.
\end{proof}

Now let's go on to dig out more about $f_{\pm1,0}(h_{2}), f_{0,\pm1}(h_{1})$. Calculate
$$\begin{aligned}d_{ 1,-1}\cdot1&=[d_{ 1,0},d_{ 0,-1}]\cdot1\\
&=d_{ 1,0}f_{ 0,-1}(h_{1})-d_{0,-1}f_{1,0}(h_{2})\\
&=(\sigma_{1}(f_{ 0,-1}(h_{1})))f_{1,0}(h_{2})-(\sigma_{2}^{-1}(f_{1,0}(h_{2})))f_{ 0,-1}(h_{1})\\
&=(\sigma_{1}-\sigma_{2}^{-1})(f_{ 0,-1}(h_{1})f_{1,0}(h_{2}));
\end{aligned}$$
$$\begin{aligned}d_{ -1,1}\cdot1&=[d_{ -1,0},d_{ 0,1}]\cdot1\\
&=d_{ -1,0}f_{ 0,1}(h_{1})-d_{0,1}f_{-1,0}(h_{2})\\
&=(\sigma_{1}^{-1}(f_{ 0,1}(h_{1})))f_{-1,0}(h_{2})-(\sigma_{2}(f_{-1,0}(h_{2})))f_{ 0,1}(h_{1})\\
&=(\sigma_{1}^{-1}-\sigma_{2})(f_{ 0,1}(h_{1})f_{-1,0}(h_{2})).
\end{aligned}$$
By $[d_{1,-1},d_{-1,1}]\cdot1=0,$ $d_{1,-1}d_{-1,1}\cdot1=d_{-1,1}d_{1,-1}\cdot1,$ we get
$$
(\sigma_{1}\sigma_{2}^{-1}((\sigma_{1}^{-1}-\sigma_{2})(f_{ 0,1}(h_{1})f_{-1,0}(h_{2}))))
((\sigma_{1}-\sigma_{2}^{-1})(f_{ 0,-1}(h_{1})f_{1,0}(h_{2})))$$
$$=(\sigma_{1}^{-1}\sigma_{2}((\sigma_{1}-\sigma_{2}^{-1})(f_{0,-1}(h_{1})f_{1,0}(h_{2}))))
((\sigma_{1}^{-1}-\sigma_{2})(f_{ 0,1}(h_{1})f_{-1,0}(h_{2}))).
$$
i.e.
$$((\sigma_{2}^{-1}-\sigma_{1})(f_{ 0,1}(h_{1})f_{-1,0}(h_{2})))
((\sigma_{1}-\sigma_{2}^{-1})(f_{ 0,-1}(h_{1})f_{1,0}(h_{2})))$$
$$=((\sigma_{2}-\sigma_{1}^{-1})(f_{0,-1}(h_{1})f_{1,0}(h_{2})))
((\sigma_{1}^{-1}-\sigma_{2})(f_{ 0,1}(h_{1})f_{-1,0}(h_{2}))).
$$
Let
$$g(h_{1},h_{2})=((\sigma_{2}-\sigma_{1}^{-1})(f_{0,-1}(h_{1})f_{1,0}(h_{2})))
((\sigma_{1}^{-1}-\sigma_{2})(f_{ 0,1}(h_{1})f_{-1,0}(h_{2}))),$$
and the above equation falls into
$$\sigma_{2}^{-1}\sigma_{1}(g(h_{1},h_{2}))=g(h_{1},h_{2}).$$

\begin{lemma} The solutions of $\sigma_{2}^{-1}\sigma_{1}(g(h_{1},h_{2}))=g(h_{1},h_{2})$ in $\mathbb{C}[h_{1},h_{2}]$ are
\begin{equation}g(h_{1},h_{2})=\mu_{0}(h_{1}+h_{2})^{m}+\mu_{1}(h_{1}+h_{2})^{m-1}+\cdot\cdot\cdot+\mu_{m-1}(h_{1}+h_{2})+\mu_{m},
\end{equation}
where $\mu_{k}\in\mathbb{C},k=0,1\cdot\cdot\cdot m$, $m\in\mathbb{Z_{+}}.$
\end{lemma}

\begin{proof} Let
$$g(h_{1},h_{2})=g(h_{1},(h_{1}+h_{2})-h_{1})=\sum_{k=0}^mb_{k}(h_{1})(h_{1}+h_{2})^{m-k}, b_{k}(h_{1})\in\mathbb{C}[h_{1}].$$
By $\sigma_{2}^{-1}\sigma_{1}(g(h_{1},h_{2}))=g(h_{1},h_{2})$, we get
$$\sum_{k=0}^mb_{k}(h_{1}-1)(h_{1}+h_{2})^{m-k}=\sum_{k=0}^mb_{k}(h_{1})(h_{1}+h_{2})^{m-k},$$
so $b_{k}(h_{1}-1)=b_{k}(h_{1})$, and $b_{k}(h_{1})$ must be constants, set $b_{k}(h_{1})=\mu_{k}$ and we get it.
\end{proof}

\begin{lemma}
$$\begin{aligned}&d_{1,0}\cdot1=\lambda_{1}(h_{2}+\alpha_{2}),d_{0,1}\cdot1=\lambda_{2}(h_{1}+\alpha_{1});\\
&d_{-1,0}\cdot1=-\lambda_{1}^{-1}(h_{2}+\alpha_{2}),d_{ 0,-1}\cdot1=-\lambda_{2}^{-1}(h_{1}+\alpha_{1}).
\end{aligned}$$
where, $\lambda_{1},\lambda_{2}\in\mathbb{C^{\ast}}$, $\alpha_{1},\alpha_{2}\in\mathbb{C}.$
\end{lemma}

\begin{proof}
By Lemma 2.3 and the previous deductions we get
\begin{equation}
\begin{aligned}
&((\sigma_{2}-\sigma_{1}^{-1})(f_{0,-1}(h_{1})f_{1,0}(h_{2})))
((\sigma_{1}^{-1}-\sigma_{2})(f_{ 0,1}(h_{1})f_{-1,0}(h_{2})))\\
&=\mu_{0}(h_{1}+h_{2})^{m}+\mu_{1}(h_{1}+h_{2})^{m-1}+\cdot\cdot\cdot+\mu_{m-1}(h_{1}+h_{2})+\mu_{m}.
\end{aligned}
\end{equation}

There must be
\begin{equation}
(\sigma_{2}-\sigma_{1}^{-1})(f_{0,-1}(h_{1})f_{1,0}(h_{2}))\in\mathbb{C}[h_{1}+h_{2}].
\end{equation}

Let
$$f_{1,0}(h_{2})=\sum_{k=0}^na_{k}(h_{1})(h_{1}+h_{2})^{n-k},$$
where $a_{0}(h_{1})\neq0$.
By (2.4), we get
$$f_{0,-1}(h_{1})\sum_{k=0}^na_{k}(h_{1})(h_{1}+h_{2}-1)^{n-k}-f_{0,-1}(h_{1}+1)\sum_{k=0}^na_{k}(h_{1})(h_{1}+h_{2})^{n-k}\in\mathbb{C}[h_{1}+h_{2}],$$
in which the leading term in $(h_{1}+h_{2})$ is
$$(f_{0,-1}(h_{1})-f_{0,-1}(h_{1}+1))a_{0}(h_{1})\in\mathbb{C}.$$
So $(f_{0,-1}(h_{1})-f_{0,-1}(h_{1}+1))\in\mathbb{C},$ and $\text {deg}f_{0,-1}(h_{1})=0,1$. Similar discussions show that
$\text {deg}f_{1,0}(h_{2})=0,1, \text {deg}f_{0,1}(h_{1})=0,1, \text {deg}f_{-1,0}(h_{2})=0,1$.

Thanks to Lemma 2.1, $d_{i,j}\cdot1\ne0$, there is only one possible case. In fact, take the following case for example,
$$\begin{aligned}&f_{1,0}(h_{2})=a_{0}, f_{ 0,1}(h_{1})=b_{0}h_{1}+b_{1}; \\
&f_{-1,0}(h_{2})=c_{0},f_{ 0,-1}(h_{1})=d_{0}h_{1}+d_{1},\end{aligned}$$
i.e.,
$$\begin{aligned}&d_{1,0}\cdot1=a_{0}, d_{0,1}\cdot1=b_{0}h_{1}+b_{1};\\
&d_{-1,0}\cdot1=c_{0},d_{ 0,-1}\cdot1=d_{0}h_{1}+d_{1},\end{aligned}$$
where, $a_{0}, b_{0}, c_{0}, d_{0}\in\mathbb{C}^{*}$ and $b_{1},d_{1}\in\mathbb{C}$.
In this case,
$$\begin{aligned}d_{1,1}\cdot1&=[d_{0,1},d_{1,0}]\cdot1=d_{0,1}a_{0}-d_{1,0}(b_{0}h_{1}+b_{1})\\
&=a_{0}(b_{0}h_{1}+b_{1})-(b_{0}h_{1}-b_{0}+b_{1})a_{0}=a_{0}b_{0}.\end{aligned}$$
i.e., $d_{1,1}\cdot1=a_{0}b_{0}$,  but then
$$d_{0,1}\cdot1=[d_{-1,0},d_{1,1}]\cdot1=d_{-1,0}a_{0}b_{0}-d_{1,1}c_{0}=a_{0}b_{0}c_{0}-c_{0}a_{0}b_{0}=0,$$
a contradiction, and this can never happen.

Each case is ruled out except
$$\text {deg}f_{1,0}(h_{2})=\text {deg}f_{0,1}(h_{1})=\text {deg}f_{-1,0}(h_{2})=\text {deg}f_{0,-1}(h_{1})=1,$$
when we set
$$\begin{aligned}&f_{1,0}(h_{2})=a_{0}h_{2}+a_{1}, f_{ 0,1}(h_{1})=b_{0}h_{1}+b_{1}; \\
&f_{-1,0}(h_{2})=c_{0}h_{2}+c_{1},f_{ 0,-1}(h_{1})=d_{0}h_{1}+d_{1},\end{aligned}$$
i.e.,
$$\begin{aligned}&d_{1,0}\cdot1=a_{0}h_{2}+a_{1}, d_{0,1}\cdot1=b_{0}h_{1}+b_{1};\\
&d_{-1,0}\cdot1=c_{0}h_{2}+c_{1},d_{ 0,-1}\cdot1=d_{0}h_{1}+d_{1}.\end{aligned}$$
Calculate
$$\begin{aligned}d_{1,1}\cdot1&=[d_{ 0,1},d_{1,0}]\cdot1\\
&=d_{ 0,1}(a_{0}h_{2}+a_{1})-d_{1,0}(b_{0}h_{1}+b_{1})\\
&=(a_{0}h_{2}+a_{1}-a_{0})(b_{0}h_{1}+b_{1})-(b_{0}h_{1}+b_{1}-b_{0})(a_{0}h_{2}+a_{1})\\
&=-a_{0}(b_{0}h_{1}+b_{1})+b_{0}(a_{0}h_{2}+a_{1})\\
&=a_{0}b_{0}(h_{2}-h_{1})+(a_{1}b_{0}-a_{0}b_{1}).
\end{aligned}$$
On the one hand,
$$\begin{aligned}d_{0,1}\cdot1&=[d_{ -1,0},d_{1,1}]\cdot1\\
&=d_{ -1,0}(a_{0}b_{0}(h_{2}-h_{1})+(a_{1}b_{0}-a_{0}b_{1}))-d_{1,1}(c_{0}h_{2}+c_{1})\\
&=(a_{0}b_{0}(h_{2}-h_{1})+(a_{1}b_{0}-a_{0}b_{1}-a_{0}b_{0}))(c_{0}h_{2}+c_{1})\\
&-(c_{0}h_{2}+c_{1}-c_{0})(a_{0}b_{0}(h_{2}-h_{1})+(a_{1}b_{0}-a_{0}b_{1}))\\
&=-a_{0}b_{0}(c_{0}h_{2}+c_{1})+c_{0}(a_{0}b_{0}(h_{2}-h_{1})+(a_{1}b_{0}-a_{0}b_{1}))\\
&=-a_{0}b_{0}c_{0}h_{1}+(a_{1}b_{0}c_{0}-a_{0}b_{1}c_{0}-a_{0}b_{0}c_{1});
\end{aligned}$$
on the other hand, $d_{0,1}\cdot1=b_{0}h_{1}+b_{1}$, so
$$ -a_{0}b_{0}c_{0}h_{1}+(a_{1}b_{0}c_{0}-a_{0}b_{1}c_{0}-a_{0}b_{0}c_{1})=b_{0}h_{1}+b_{1},$$
and we get
\begin{equation}
  \left\{
   \begin{aligned}
   a_{0}c_{0}=-1,  \\
   \frac{a_{1}}{a_{0}}= \frac{c_{1}}{c_{0}}.\\
   \end{aligned}
   \right.
  \end{equation}\

Similarly, consider on the one hand,
$$d_{1,0}\cdot1=[d_{1,1},d_{0,-1}]\cdot1=-a_{0}b_{0}d_{0}h_{2}+(a_{0}b_{1}d_{0}-a_{1}b_{0}d_{0}-a_{0}b_{0}d_{1});$$
on the other hand,
$$d_{1,0}\cdot1=a_{0}h_{2}+a_{1},$$
so
$$-a_{0}b_{0}d_{0}h_{2}+(a_{0}b_{1}d_{0}-a_{1}b_{0}d_{0}-a_{0}b_{0}d_{1})=a_{0}h_{2}+a_{1},$$
and we get
\begin{equation}
  \left\{
   \begin{aligned}
   b_{0}d_{0}=-1,  \\
   \frac{b_{1}}{b_{0}}= \frac{d_{1}}{d_{0}}.\\
   \end{aligned}
   \right.
  \end{equation}

By (2.5) and (2.6), let
$$\begin{aligned}&a_{0}=\lambda_{1},a_{1}=\alpha_{2}\lambda_{1},c_{0}=-\lambda_{1}^{-1},c_{1}=-\alpha_{2}\lambda_{1}^{-1};\\
&b_{0}=\lambda_{2},b_{1}=\alpha_{1}\lambda_{2},d_{0}=-\lambda_{2}^{-1},d_{1}=-\alpha_{1}\lambda_{2}^{-1},\end{aligned}$$
and the lemma holds.
\end{proof}

\begin{lemma}For $(i,j)\in\mathbb{Z}^{2}\setminus\{(0,0)\}$,
\begin{equation}d_{i,j}\cdot1=(-1)^{i+1}\lambda_{1}^{i}\lambda_{2}^{j}(i(h_{2}+\alpha_{2})-j(h_{1}+\alpha_{1})), (i,j)\neq(0,0).\end{equation}
\end{lemma}

\begin{proof} This follows from Lemma 2.4 and Definition 1.

\end{proof}

Consequently we have

\begin{theorem} Let $M$ be a $\tilde{S}_{2}$-module which is a free  $U(\h_2)$-module of rank 1, on which $\tilde{S}_{2}'$ acts nontrivially.
Then $M\cong \Omega(\lambda_{1},\lambda_{2},\alpha_{1},\alpha_{2})$ for some $\lambda_{1}, \lambda_{2}\in \C^*, \alpha_{1}, \alpha_{2}\in \C$.
\end{theorem}

Set $$\Omega'(\lambda_{1},\lambda_{2},\alpha_{1},\alpha_{2})=\langle h_{1}+\alpha_{1}, h_{2}+\alpha_{2}\rangle,$$
the ideal of the polynomial algebra $\C[h_{1},h_{1}]$ generated by $h_{1}+\alpha_{1}$ and $ h_{2}+\alpha_{2}$. $\Omega'(\lambda_{1},\lambda_{2},\alpha_{1},\alpha_{2})$ is a submodule considering (2.7). For the simplicity, we provide

\begin{theorem}$\Omega'(\lambda_{1},\lambda_{2},\alpha_{1},\alpha_{2})$ is the unique irreducible $\tilde{S}_{2}$-submodule of $\Omega(\lambda_{1},\lambda_{2},\alpha_{1},\alpha_{2})$, and
$$\Omega(\lambda_{1},\lambda_{2},\alpha_{1},\alpha_{2})/\Omega'(\lambda_{1},\lambda_{2},\alpha_{1},\alpha_{2})\cong\mathbb{C}.$$
\end{theorem}

\begin{proof} The isomorphism is immediate. For the irreducibility, take any nonzero $F(h_{1},h_{2})\in \Omega(\lambda_{1},\lambda_{2},\alpha_{1},\alpha_{2})$, let $N=U(\tilde{S}_{2})(F(h_{1},h_{2}))$, it's sufficient to prove $h_{1}+\alpha_{1},h_{2}+\alpha_{2}\in N$. By looking at
\begin{equation}
d_{i,0}F(h_{1},h_{2})=(-1)^{i+1}i\lambda_{1}^{i}F(h_{1}-i,h_{2})(h_{2}+\alpha_{2}), i\in \mathbb{Z},
\end{equation}
we see that there is a nonzero $G(h_{2})\in\mathbb{C}[h_{2}]\cap N$. Considering
$$d_{0,j}G(h_{2})=j\lambda_{2}^{j}G(h_{2}-j)(h_{1}+\alpha_{1}), j\in \mathbb{Z},$$
we see that $h_{1}+\alpha_{1}\in N$. By replacing $F(h_{1},h_{2})$ with $h_{1}+\alpha_{1}$ in (2.8) we will get $h_{2}+\alpha_{2}\in N$.
In fact, $\Omega'(\lambda_{1},\lambda_{2},\alpha_{1},\alpha_{2})$ is the unique nontrivial submodule, the uniqueness is clear by the proof of the irreducibility.
\end{proof}

The simple module $\Omega'(\lambda_{1},\lambda_{2},\alpha_{1},\alpha_{2})$ is completely determined by the quadruple $(\lambda_{1},\lambda_{2},\alpha_{1},\alpha_{2})$.
\begin{theorem}  Let $\lambda_{1}, \lambda_{2}, \lambda_{1}', \lambda_{2}'\in \C^*, \alpha_{1}, \alpha_{2}, \alpha_{1}', \alpha_{2}'\in \C$.
Then the simple $\tilde{S}_{2}$-modules $\Omega'(\lambda_{1},\lambda_{2},\alpha_{1},\alpha_{2})\cong\Omega'(\lambda_{1}',\lambda_{2}',\alpha_{1}',\alpha_{2}')$ if and only if $\lambda_{1}=\lambda_{1}',\lambda_{2}=\lambda_{2}',\alpha_{1}=\alpha_{1}', \alpha_{2}=\alpha_{2}'$.
\end{theorem}
\begin{proof} The sufficiency is evident. For the necessity, let $$\psi:\Omega'(\lambda_{1},\lambda_{2},\alpha_{1},\alpha_{2})\rightarrow\Omega'(\lambda_{1}',\lambda_{2}',\alpha_{1}',\alpha_{2}')$$
be an isomorphism between the two $\tilde{S}_{2}$-modules. Let's consider
$$(d_{0,1}-\lambda_{2}(h_{1}+\alpha_{1}))(h_{1}+\alpha_{1})=0,$$
letting $\psi$ act on it and we get
$$\begin{aligned}
&(d_{0,1}-\lambda_{2}(h_{1}+\alpha_{1}))\psi(h_{1}+\alpha_{1})=0,\\
&\lambda_{2}'(\sigma_{2}(\psi(h_{1}+\alpha_{1})))(h_{1}+\alpha_{1}')=\lambda_{2}(h_{1}+\alpha_{1})\psi(h_{1}+\alpha_{1}),
\end{aligned}$$
so $\lambda_{2}=\lambda_{2}'$. Rewrite the previous equality
$$(\sigma_{2}(\psi(h_{1}+\alpha_{1}))-\psi(h_{1}+\alpha_{1}))(h_{1}+\alpha_{1}')=(\alpha_{1}-\alpha_{1}')\psi(h_{1}+\alpha_{1}),$$
we see that $\psi(h_{1}+\alpha_{1})\in\mathbb{C}[h_{1}]$ and $\alpha_{1}=\alpha_{1}'$. By similar discussions, we can show that
$\lambda_{1}=\lambda_{1}', \alpha_{2}=\alpha_{2}'$ and $\psi(h_{2}+\alpha_{2})\in\mathbb{C}[h_{2}]$.
\end{proof}

We remark that the $\tilde{S}_{2}$-modules restricted from the $\W_2$-modules $\Omega(\Lambda_2,a)$ stated in Theorem 9([TZ2]) exhaust all the $\mathscr{L}$-modules $\Omega(\lambda_{1},\lambda_{2},0,0)$,  where $\Lambda_2=(\lambda_{1},\lambda_{2})\in (\C^*)^2, a\in\C$.

\subsection{$\tilde{S}_{n}$-module structures  on $U(\h_n)(n\geq3)$}

\vskip 5pt

First, let's collect some basic facts for $\tilde{S}_{n}$.

{\bf{Fact 1}}.  $\tilde{S}_{n}=\text {Span}_{\mathbb{C}}\{\partial_i,r_{j}t^{r}\partial_i-r_{i}t^{r}\partial_j\mid i,j\in\{1,\cdot\cdot\cdot,n\},r\in \Z^n\}$ .

{\bf{Fact 2}}.  $\tilde{S}_{n}=\tilde{S}_{n}'\oplus\mathfrak{h}_{n}$,
$$\tilde{S}_{n}'=\text {Span}_{\mathbb{C}}\{r_{j}t^{r}\partial_i-r_{i}t^{r}\partial_j\mid i,j\in\{1,\cdot\cdot\cdot,n\},r\in \Z^n\}.$$

{\bf{Fact 3}}.  $\tilde{S}_{n}$ is generated by $\{\partial_i,t_{i}\partial_j,t_{i}^{-1}\partial_j\mid i,j\in\{1,\cdot\cdot\cdot,n\}, i\neq j\}$.

{\bf{Fact 4}}.  $\tilde{S}_{n}'$ is generated by $\{t_{i}\partial_j,t_{i}^{-1}\partial_j\mid i,j\in\{1,\cdot\cdot\cdot,n\}, i\neq j\}$.

It is easy to see that {\bf{Fact 3}} follows from {\bf{Fact 4}}, while {\bf{Fact 4}} can be established in a similar way as that of Lemma 3.11.

Let $\tilde{S}_{ij}$ be the subalgebra of $\tilde{S}_{n}$ generated by
$$\{\partial_i,\partial_j,t_{i}\partial_j,t_{j}\partial_i,t_{i}^{-1}\partial_j,t_{j}^{-1}\partial_i\mid i,j\in\{1,\cdot\cdot\cdot,n\}, i\neq j\}.$$
Note that $\tilde{S}_{ij}=\tilde{S}_{ji}$, there are altogether
$\left(
\begin{array}{c}
n\\
2\\
\end{array}
\right)$
such subalgebras in $\tilde{S}_{n}$.

{\bf{Fact 5}}. $\tilde{S}_{ij}\cong\mathscr{L}$, with the isomorphism:
$$\begin{aligned}
\phi_{ij}:\tilde{\mathrm{S}}_{ij}&\rightarrow\mathscr{L}\\
\partial_i&\mapsto h_{1}\\
\partial_j&\mapsto h_{2}\\
t_{i}\partial_j&\mapsto d_{1,0} \\
-t_{i}^{-1}\partial_j&\mapsto d_{-1,0}\\
t_{j}\partial_i&\mapsto d_{0, 1} \\
-t_{j}^{-1}\partial_i&\mapsto d_{0, -1}\\
t_{i}^{r_{i}}\partial_j&\mapsto \frac{(-1)^{r_{i}+1}}{r_{i}}d_{r_{i},0}\\
t_{j}^{r_{j}}\partial_i&\mapsto \frac{1}{r_{j}}d_{0, r_{j}}\\
r_{j}t_{i}^{r_{i}}t_{j}^{r_{j}}\partial_i-r_{i}t_{i}^{r_{i}}t_{j}^{r_{j}}\partial_j&\mapsto(-1)^{r_{i}}d_{r_{i},r_{j}},\\
\end{aligned}$$
where $(r_{i}, r_{j})\in \Z^2$.

{\bf{Fact 6}}. $\tilde{S}_{ij}=\tilde{S}_{ij}'\oplus\mathbb{C}\partial_i\oplus\mathbb{C}\partial_j$,
$\tilde{S}_{ij}'\cong\mathscr{L}'$, and $\tilde{S}_{ij}'$ is simple.

\begin{lemma}
$\tilde{S}_{n}'=\left(\text {ad}\tilde{S}_{n}'\right)\left(\tilde{S}_{ij}'\right)$ and $\tilde{S}_{n}'$ is simple.
\end{lemma}

\begin{proof}
$\tilde{S}_{n}'=\left(\text {ad}\tilde{S}_{n}'\right)\left(\tilde{S}_{ij}'\right)$ due to
$$r_{j}t^{r}\partial_i-r_{i}t^{r}\partial_j
=[t_{1}^{r_{1}}\cdot\cdot\cdot \overset{\wedge}{t_{j}^{r_{j}}}\cdot\cdot\cdot t_{n}^{r_{n}}\partial_j,t_{j}^{r_{j}}\partial_i].$$
For the simplicity,
let $I$ be a nonzero ideal of $\tilde{S}_{n}'$, by {\bf{Fact 4}},  $\{\tilde{S}_{ij}'\mid i,j\in\{1,\cdot\cdot\cdot,n\}, i\neq j\}$ generate $\tilde{S}_{n}'$, then there must be some  $\tilde{S}_{ij}'$  such that $I\cap\tilde{S}_{ij}'=\tilde{S}_{ij}'$, then $I=\tilde{S}_{n}'$.
\end{proof}

Let $M=U(\h_n)1\cong\mathbb{C}[\partial_1,\partial_2,\cdot\cdot\cdot,\partial_n]$ $(3\le n \le \infty)$ be a $\tilde{S}_{n}$-module
with left multiplication action of $\h_n$, on which $\tilde{S}_{n}'$ acts nontrivially. To get a clear picture of it, let's do some preparations.
\begin{lemma} For any $f(\partial_1,\partial_2,...,\partial_n)\in M$, $r=(r_1,r_2,\cdot\cdot\cdot,r_n)\in \mathbb{Z}^{n}$,
$i,j\in\{1,\cdot\cdot\cdot,n\}, i\neq j$, we have
$$\begin{aligned}
&(r_{j}t^{r}\partial_i-r_{i}t^{r}\partial_j)\cdot f(\partial_1,\partial_2,...,\partial_n)1\\
&=f(\partial_1-r_1,\partial_2-r_2,...,\partial_n-r_n)((r_{j}t^{r}\partial_i-r_{i}t^{r}\partial_j)\cdot 1),
\end{aligned}$$
 and $(r_{j}t^{r}\partial_i-r_{i}t^{r}\partial_j)\cdot 1\ne0$.
\end{lemma}

\begin{proof} The proof is similar to that of Lemma 2.1.
\end{proof}

For  $i,j\in\{1,\cdot\cdot\cdot,n\}, i\neq j$, set
$$R_{ij}=\mathbb{C}[\partial_1,\cdot\cdot\cdot,\partial_{i-1},\partial_{i+1},\cdot\cdot\cdot,\partial_{j-1},\partial_{j+1},\cdot\cdot\cdot,\partial_n].$$

\begin{lemma} For $i,j\in\{1,\cdot\cdot\cdot,n\}, i\neq j$,
there exist $\lambda_{i}\in\mathbb{C}^{\ast}, \alpha_{j}\in\mathbb{C}$,such that
\begin{equation}
t_{i}^{r_{i}}\partial_j\cdot1=\lambda_{i}^{r_{i}}(\partial_j+\alpha_{j}), \forall r_{i}\in\mathbb{Z}^{\ast}.
\end{equation}
\end{lemma}

\begin{proof} Note that $n\geq3$. By the isomorphism $\phi_{ij}$ in {\bf{Fact 5}} and (2.7), with $h_{2}$ replaced by $\partial_j$, $\alpha_{2}$ by $\alpha_{ij}$ and $\lambda_{1}$ by $\lambda_{ij}$, we get
$$t_{i}^{r_{i}}\partial_j\cdot1= \frac{(-1)^{r_{i}+1}}{r_{i}}d_{r_{i},0}\cdot1=\lambda_{ij}^{r_{i}}(\partial_j+\alpha_{ij}),\alpha_{ij}\in R_{ij}.$$

{\bf{Claim 1}}      $\alpha_{ij}\in\mathbb{C}$.

For arbitrary pairwise different $i,j,k\in\{1,\cdot\cdot\cdot,n\}$, set
$$t_{i}\partial_j\cdot1=\lambda_{ij}(\partial_j+\alpha_{ij}), t_{k}\partial_j\cdot1=\lambda_{kj}(\partial_j+\alpha_{kj}),$$
where $\alpha_{ij}\in R_{ij}, \alpha_{kj}\in R_{kj}$.
Consider
$$\begin{aligned}
&[t_{k}\partial_j, t_{i}\partial_j]\cdot1=0,\\
&t_{k}\partial_jt_{i}\partial_j\cdot1=t_{i}\partial_jt_{k}\partial_j\cdot1,\\
&\lambda_{ij}\lambda_{kj}(\partial_j+\sigma_{k}(\alpha_{ij}))(\partial_j+\alpha_{kj})=\lambda_{ij}\lambda_{kj}(\partial_j+\alpha_{ij})(\partial_j+\sigma_{i}(\alpha_{kj})),\\
&(\partial_j+\sigma_{k}(\alpha_{ij}))(\partial_j+\alpha_{kj})=(\partial_j+\alpha_{ij})(\partial_j+\sigma_{i}(\alpha_{kj})),\\
&(\sigma_{k}(\alpha_{ij})+\alpha_{kj})\partial_j+\sigma_{k}(\alpha_{ij})\alpha_{kj}=(\sigma_{i}(\alpha_{kj})+\alpha_{ij})\partial_j+\sigma_{i}(\alpha_{kj})\alpha_{ij},
\end{aligned}$$
from the last equality, we deduce that
$$\left\{
\begin{aligned}
\sigma_{k}(\alpha_{ij})-\alpha_{ij} & =  \sigma_{i}(\alpha_{kj})- \alpha_{kj}\\
\sigma_{k}(\alpha_{ij})\alpha_{kj} & =  \sigma_{i}(\alpha_{kj})\alpha_{ij}
\end{aligned}
\right.
,$$
so $(\sigma_{k}(\alpha_{ij})-\alpha_{ij})\alpha_{kj}=( \sigma_{i}(\alpha_{kj})-\alpha_{kj})\alpha_{ij}$. If $\sigma_{k}(\alpha_{ij})-\alpha_{ij}\neq0$ for some $k$, then $\alpha_{kj}=\alpha_{ij}\in R_{ij}\cap R_{kj}$, and $\sigma_{k}(\alpha_{ij})-\alpha_{ij}=0$, a  contradiction. Thus $\sigma_{k}(\alpha_{ij})-\alpha_{ij}=0, \forall k\in\{1,\cdot\cdot\cdot,n\}, k\neq i,j$, and $\alpha_{ij}\in\C$.

{\bf{Claim 2}} The $\widetilde{\mathrm{S}}_{ij}$-actions on $M$ are compatible.
For
$$\begin{aligned}
&t_{i}\partial_j\cdot1=\lambda_{ij}(\partial_j+\alpha_{ij}), t_{k}\partial_j\cdot1=\lambda_{kj}(\partial_j+\alpha_{kj});\\
&t_{i}\partial_k\cdot1=\lambda_{ik}(\partial_k+\alpha_{ik}),t_{k}\partial_i\cdot1=\lambda_{ki}(\partial_i+\alpha_{ki}),
\end{aligned}$$
calculate
$$\begin{aligned}
t_{i}t_{k}\partial_j\cdot1&=[t_{k}\partial_i, t_{i}\partial_j]\cdot1\\
&=t_{k}\partial_i t_{i}\partial_j\cdot1-t_{i}\partial_j t_{k}\partial_i\cdot1\\
&=\lambda_{ij}t_{k}\partial_i(\partial_j+\alpha_{ij})-\lambda_{ki}t_{i}\partial_j(\partial_i+\alpha_{ki})\\
&=\lambda_{ij}\lambda_{ki}(\partial_j+\alpha_{ij})(\partial_i+\alpha_{ki})-\lambda_{ki}\lambda_{ij}(\partial_i+\alpha_{ki}-1)(\partial_j+\alpha_{ij})\\
&=\lambda_{ij}\lambda_{ki}(\partial_j+\alpha_{ij}),
\end{aligned}$$

$$\begin{aligned}
t_{i}t_{k}\partial_j\cdot1&=[t_{i}\partial_k, t_{k}\partial_j]\cdot1\\
&=t_{i}\partial_k t_{k}\partial_j\cdot1-t_{k}\partial_j t_{i}\partial_k\cdot1\\
&=\lambda_{kj}t_{i}\partial_k(\partial_j+\alpha_{kj})-\lambda_{ik}t_{k}\partial_j(\partial_k+\alpha_{ik})\\
&=\lambda_{ik}\lambda_{kj}(\partial_j+\alpha_{kj})(\partial_k+\alpha_{ik})-\lambda_{ik}\lambda_{kj}(\partial_k+\alpha_{ik}-1)(\partial_j+\alpha_{kj})\\
&=\lambda_{ik}\lambda_{kj}(\partial_j+\alpha_{kj}),
\end{aligned}$$
there must be $\alpha_{ij}=\alpha_{kj}$, which will be denoted by $\alpha_{j}$. Besides we also see that $\lambda_{ij}\lambda_{ki}=\lambda_{ik}\lambda_{kj},\frac{\lambda_{ik}}{\lambda_{ij}}=\frac{\lambda_{ki}}{\lambda_{kj}}$ .
To go on, let
$$t_{j}\partial_i\cdot1=\lambda_{ji}(\partial_i+\alpha_{i}), t_{j}\partial_k\cdot1=\lambda_{jk}(\partial_k+\alpha_{k}),$$
similar discussions will show that
$$\frac{\lambda_{jk}}{\lambda_{ji}}=\frac{\lambda_{ik}}{\lambda_{ij}}, \frac{\lambda_{jk}}{\lambda_{ji}}=\frac{\lambda_{kj}}{\lambda_{ki}},$$
so $\frac{\lambda_{ik}}{\lambda_{ij}}=\frac{\lambda_{ki}}{\lambda_{kj}}=\frac{\lambda_{jk}}{\lambda_{ji}}=\pm1$. We claim that $\frac{\lambda_{ik}}{\lambda_{ij}}=\frac{\lambda_{ki}}{\lambda_{kj}}=\frac{\lambda_{jk}}{\lambda_{ji}}=1$, for which investigate into

$$\begin{aligned}
(t_{i}t_{j}t_{k}\partial_i-t_{i}t_{j}t_{k}\partial_k)\cdot1&=[t_{i}t_{j}\partial_k,t_{k}\partial_i]\cdot1\\
&=(t_{i}t_{j}\partial_k)t_{k}\partial_i\cdot1-t_{k}\partial_i(t_{i}t_{j}\partial_k)\cdot1\\
&=(t_{i}t_{j}\partial_k)\lambda_{ki}(\partial_i+\alpha_{i})-t_{k}\partial_i \lambda_{ij}\lambda_{jk}(\partial_k+\alpha_{k})\\
&=\lambda_{ij}\lambda_{jk}\lambda_{ki}(\partial_k+\alpha_{k})(\partial_i+\alpha_{i}-1)\\
&-\lambda_{ij}\lambda_{jk}\lambda_{ki}(\partial_k+\alpha_{k}-1)(\partial_i+\alpha_{i})\\
&=\lambda_{ij}\lambda_{jk}\lambda_{ki}((\partial_i+\alpha_{i})-(\partial_k+\alpha_{k}));
\end{aligned}$$
and
$$\begin{aligned}
(t_{i}t_{j}t_{k}\partial_i-t_{i}t_{j}t_{k}\partial_k)\cdot1&=[t_{i}\partial_k,t_{j}t_{k}\partial_i]\cdot1\\
&=t_{i}\partial_k(t_{j}t_{k}\partial_i)\cdot1-(t_{j}t_{k}\partial_i)t_{i}\partial_k\cdot1\\
&=t_{i}\partial_k \lambda_{ji}\lambda_{kj}(\partial_i+\alpha_{i})-(t_{j}t_{k}\partial_i)\lambda_{ik}(\partial_k+\alpha_{k})\\
&=\lambda_{ik}\lambda_{ji}\lambda_{kj}(\partial_k+\alpha_{k})(\partial_i+\alpha_{i}-1)\\
&-\lambda_{ik}\lambda_{ji}\lambda_{kj}(\partial_k+\alpha_{k}-1)(\partial_i+\alpha_{i})\\
&=\lambda_{ik}\lambda_{ji}\lambda_{kj}((\partial_i+\alpha_{i})-(\partial_k+\alpha_{k})),
\end{aligned}$$
we get $\lambda_{ij}\lambda_{jk}\lambda_{ki}=\lambda_{ik}\lambda_{ji}\lambda_{kj}$, so $\frac{\lambda_{ik}}{\lambda_{ij}}=\frac{\lambda_{ki}}{\lambda_{kj}}=\frac{\lambda_{jk}}{\lambda_{ji}}=1$, and $\lambda_{ik}=\lambda_{ij}$, which will be denoted by $\lambda_{i}$.
\end{proof}
\begin{lemma} For  $i,j\in\{1,\cdot\cdot\cdot,n\}, i\neq j, r=(r_{1}, \cdot\cdot\cdot, r_{n})\in \Z^n$,
$$\begin{aligned}
&(r_{j}t^{r}\partial_i-r_{i}t^{r}\partial_j)f(\partial_1,\partial_2,\cdot\cdot\cdot\partial_n)\\
&=f(\partial_1-r_{1},\partial_2-r_{2},\cdot\cdot\cdot\partial_n-r_{n})((r_{j}t^{r}\partial_i-r_{i}t^{r}\partial_j)\cdot1)\\
&=\lambda_{1}^{r_{1}}\cdot\cdot\cdot\lambda_{j}^{r_{j}}\cdot\cdot\cdot \lambda_{n}^{r_{n}}(r_{j}(\partial_i+\alpha_{i})-r_{i}(\partial_j+\alpha_{j}))
f(\partial_1-r_{1},\partial_2-r_{2},\cdot\cdot\cdot\partial_n-r_{n}).
\end{aligned}$$
\end{lemma}
\begin{proof}
The first equation is due to Lemma 2.10, while the second is true by Lemma 2.11 and Definition 1.
\end{proof}

Consequently,

\begin{theorem} Let $M$ be a $\tilde{S}_{n}$-module which is a free $U(\h_n)$-module of rank 1, on which $\tilde{S}_{n}'$ acts nontrivially.
Then $M\cong\Omega(\Lambda_n,\alpha)$ for some $\Lambda_n\in (\C^*)^n, \alpha\in (\C)^{n}$.
\end{theorem}

Set $$\Omega'(\Lambda_n,\alpha)=\langle \partial_1+\alpha_{1}, \cdot\cdot\cdot, \partial_n+\alpha_{n}\rangle,$$
the ideal of the polynomial algebra $\C[\partial_1, \cdot\cdot\cdot, \partial_n]$ generated by $\{\partial_i+\alpha_{i}\mid i=1, \cdot\cdot\cdot, n\}$. $\Omega'(\Lambda_n,\alpha)$ is a submodule of $\Omega(\Lambda_n,\alpha)$. By similar discussions, we can prove that
\begin{theorem}
$\Omega(\Lambda_n,\alpha)$ has a unique irreducible submodule $\Omega'(\Lambda_n,\alpha)$, and the quotient is $\mathbb{C}$.
\end{theorem}

For  $i\in\{1,\cdot\cdot\cdot,n\}$, set
$$R_{i}=\mathbb{C}[\partial_1,\cdot\cdot\cdot,\partial_{i-1},\partial_{i+1},\cdot\cdot\cdot,\partial_n].$$

\begin{theorem}Let $\Lambda_n, \Lambda_n'\in (\C^*)^n, \alpha, \alpha\in (\C)^{n}$. Then
the simple $\tilde{S}_{n}$-modules $\Omega'(\Lambda_n,\alpha)\cong \Omega'(\Lambda'_n,\alpha')$ if and only if $\Lambda_n=\Lambda'_n,\alpha=\alpha'$.
\end{theorem}
\begin{proof} The sufficiency is evident. For the necessity, let $\psi:\Omega'(\Lambda_n,\alpha))\rightarrow \Omega'(\Lambda'_n,\alpha')$
be an isomorphism between the two $\tilde{S}_{n}$-modules. Let's consider
$$(t_{j}\partial_{i}-\lambda_{j}(\partial_{i}+\alpha_{i}))(\partial_{i}+\alpha_{i})=0, i,j\in\{1,\cdot\cdot\cdot,n\}, i\neq j,$$
letting $\psi$ act on it and we get
$$\begin{aligned}
&(t_{j}\partial_{i}-\lambda_{j}(\partial_{i}+\alpha_{i}))\psi(\partial_{i}+\alpha_{i})=0,\\
&\lambda_{j}'(\sigma_{j}(\psi(\partial_{i}+\alpha_{i})))(\partial_{i}+\alpha_{i}')=\lambda_{j}(\partial_{i}+\alpha_{i})\psi(\partial_{i}+\alpha_{i}),
\end{aligned}$$
so $\lambda_{j}=\lambda_{j}'$. Rewrite the previous equality
$$(\sigma_{j}(\psi(\partial_{i}+\alpha_{i}))-\psi(\partial_{i}+\alpha_{i}))(\partial_{i}+\alpha_{i}')=(\alpha_{i}-\alpha_{i}')\psi(\partial_{i}+\alpha_{i}),$$
we see that $\psi(\partial_{i}+\alpha_{i})\in R_{j}$ and $\alpha_{i}=\alpha_{i}'$. Furthermore, we see that
$\psi(\partial_{j}+\alpha_{j})\in R_{i}$.
\end{proof}
We also remark that when $\alpha_{i}=0, i=1,\cdot\cdot\cdot,n$, the $\tilde{S}_{n}$-modules $\Omega(\Lambda_n,\mathbf{0})$ is the restriction of $\Omega(\Lambda_n,a)$ ([TZ1]) from $\W_n$.\\

\section{$\bar{S}_{n}$-module structures  on $U(\h_n)$}

\subsection{$\bar{S}_{2}$-module structures  on $U(\h_2)$}
First, let's work on the $\bar{S}_{2}$-module structures  on $U(\h_2)$.
It's convenient to quote the following  notations:
$$\begin{aligned}
&h_{1}=\partial_1, h_{2}=\partial_2, \h_{2}=\mathbb{C}h_1\oplus\mathbb{C}h_2;\\
&l_{r}=l_{r_{1},r_{2}}=(r_{2}+1)t^{r}\partial_1-(r_{1}+1)t^{r}\partial_2,
\end{aligned}$$
then
$$\bar{S}_{2}=\text{Span}_{\mathbb{C}}\{h_{1}, h_{2}, l_{r}\mid r=(r_{1},r_{2})\in \mathbb{Z}^{2}, r\neq(-1,-1),(-2,-2)\},$$
with a Cartan subalgebra $\h_{2}=\mathbb{C}h_1\oplus\mathbb{C}h_2$, and the derived algebra
$$\bar{S}_{2}'=\text{Span}_{\mathbb{C}}\{l_{r}\mid r=(r_{1},r_{2})\in \mathbb{Z}^{2}, r\neq(-1,-1),(-2,-2)\},$$
whose brackets looks like
$$\begin{aligned}
&[h_{i},l_{r}]=r_{i}l_{r}, i=1,2;\\
&[l_{r},l_{r'}]=\left|
 \begin{array}{cc}
  r'_{1}+1 &  r_{1}+1 \\
  r'_{2}+1 &  r_{2}+1
  \end{array}
 \right|l_{r+r'},
\end{aligned}$$
where $r=(r_{1},r_{2}), r'=(r'_{1},r'_{2})\in \mathbb{Z}^{2}$ with $r, r'\neq(-1,-1),(-2,-2)$, and $l_{0,0}=h_{1}-h_{2}$, while $l_{-1,-1}=0$ understandably.

Let's consider the $\bar{S}_{2}$- module structures on $M=U(\h_{2})\cdot1\cong\mathbb{C}[h_{1},h_{2}]$ that is a free $U(\h_{2})$-module of rank 1, on which $\bar{S}_{2}'$ acts nontrivially.

Note that $\overline{S}_{2}'\cap \h_{2}=\mathbb{C}(h_{1}-h_2)$  and $\overline{S}_{2}'$ is simple, so a free $U(\h_{2})$-module of rank 1 is naturally a module on which $\overline{S}_{2}'$ acts nontrivially.

\begin{lemma} For any $f(h_{1},h_{2})\cdot1\in M$, where $r=(r_{1},r_{2})\in\mathbb{Z}^{2}$ with $r\neq(-1,-1),(-2,-2)$ we have
\begin{equation}
l_{r}\cdot f(h_{1},h_{2})1=f(h_{1}-r_{1},h_{2}-r_{2})(l_{r}\cdot1),
\end{equation}
and $l_{r}\cdot1\ne0$.
\end{lemma}
\begin{proof} The proof is similar to that of Lemma 2.1.
\end{proof}\

Straightforward calculations show that $\bar{S}_{2}'$ is generated by $l_{\pm1,0}, l_{0,\pm1}, l_{-2,0},l_{0,-2}.$
By Lemma 3.1, the $\bar{S}_{2}-$module structures on $M$ is determined by
$$l_{\pm1,0}\cdot1, l_{0,\pm1}\cdot1, l_{-2,0}\cdot1, l_{0,-2}\cdot1.$$
There are polynomials
$$f_{\pm1,0}(h_{1},h_{2}), f_{0,\pm1}(h_{1},h_{2}), f_{-2,0}(h_{1},h_{2}),f_{0-2}(h_{1},h_{2})\in\mathbb{C}[h_{1},h_{2}]$$
such that
$$\begin{aligned}
&l_{\pm1,0}\cdot1=f_{\pm1,0}(h_{1},h_{2}), l_{0,\pm1}\cdot1=f_{0,\pm1}(h_{1},h_{2}), \\
&l_{-2,0}\cdot1=f_{-2,0}(h_{1},h_{2}),l_{0,-2}\cdot1=f_{0-2}(h_{1},h_{2}).
\end{aligned}$$
For further convenience, denote $l_{r_{1},r_{2}}\cdot1=f_{r_{1},r_{2}}(h_{1},h_{2})=f_{r_{1},r_{2}}\in\mathbb{C}[h_{1},h_{2}]$.

Calculate
$$\begin{aligned}
&[l_{1,0},l_{-1,0}]=\left|
                    \begin{array}{cc}
                      0 & 2 \\
                      1 & 1
                    \end{array}
                    \right|l_{0,0}=-2l_{0,0},\\
&[l_{1,0},l_{-1,0}]\cdot1=-2(h_{1}-h_{2}),\\
&l_{1,0}f_{-1,0}(h_{1},h_{2})-l_{-1,0}f_{1,0}(h_{1},h_{2})=-2(h_{1}-h_{2}),\\
&(\sigma_{1}f_{-1,0})f_{1,0}-(\sigma_{1}^{-1}f_{1,0})f_{-1,0}=-2(h_{1}-h_{2}),\\
&\sigma_{1}((\sigma_{1}^{-1}f_{1,0})f_{-1,0})-(\sigma_{1}^{-1}f_{1,0})f_{-1,0}=-2(h_{1}-h_{2}).
\end{aligned}$$

Set
$$(\sigma_{1}^{-1}f_{1,0})f_{-1,0}=\sum_{k=0}^{m} a_{k}(h_{2})(h_{1}-h_{2})^{k},$$
we get
$$\sum_{k=0}^{m} a_{k}(h_{2})((h_{1}-h_{2})^{k}-(h_{1}-h_{2}-1)^{k})=2(h_{1}-h_{2}),$$
and
$$a_{2}(h_{2})((h_{1}-h_{2})^{2}-(h_{1}-h_{2}-1)^{2})+a_{1}(h_{2})((h_{1}-h_{2})-(h_{1}-h_{2}-1))=2(h_{1}-h_{2});$$
compare the coefficients,
$$\left\{
\begin{aligned}
2a_{2}(h_{2}) & = & 2 \\
a_{1}(h_{2})-a_{2}(h_{2})& = & 0
\end{aligned}
\right.
,$$
so
$$\left\{
\begin{aligned}
a_{2}(h_{2}) & = & 1 \\
a_{1}(h_{2})& = & 1
\end{aligned}
\right.
$$
and
\begin{equation}
(\sigma_{1}^{-1}f_{1,0})f_{-1,0}=(h_{1}-h_{2})^{2}+(h_{1}-h_{2})+a_{0}(h_{2}).
\end{equation}
Thanks to (3.2), we have $\text {deg}_{1}f_{\pm1,0}\leq 2$. By similar discussions, we can prove $\text {deg}_{2}f_{0,\pm1}\leq 2$.
Furthermore, we can get the following lemma:
\begin{lemma}
$\text {deg}_{1}f_{\pm1,0}=1,\text {deg}_{2}f_{0,\pm1}=1.$
\end{lemma}

\begin{proof}
If $f_{1,0}=a\in\C^*$, then
$$\begin{aligned}
f_{n+1,-1}=l_{n+1,-1}\cdot1
&=\frac{1}{n+1}[l_{1,0},l_{n,-1}]\cdot1\\
&=\frac{1}{n+1}((\sigma_{1}f_{n,-1})f_{1,0}-(\sigma_{1}^{n}\sigma_{2}^{-1}f_{1,0})f_{n,-1})\\
&=\frac{a}{n+1}(\sigma_{1}f_{n,-1}-f_{n,-1}),
\end{aligned}$$
note that  $\text {deg}_{1}f_{n+1,-1}=\text {deg}_{1}f_{n,-1}-1$, so there must be $f_{n,1}=0$ for some $n\in\mathbb{N}$, a contradiction.

If $f_{-1,0}=a\in\C^*$, then
$$\begin{aligned}
f_{-(n+1),0}=l_{-(n+1),0}\cdot1
&=\frac{1}{1-n}[l_{-1,0},l_{-n,0}]\cdot1\\
&=\frac{1}{1-n}((\sigma_{1}^{-1}f_{-n,0})f_{-1,0}-(\sigma_{1}^{-n}f_{-1,0})f_{-n,0})\\
&=\frac{a}{1-n}(\sigma_{1}^{-1}f_{-n,0}-f_{-n,0}),
\end{aligned}$$
note that  $\text {deg}_{1}f_{-(n+1),0}=\text {deg}_{1}f_{-n,0}-1$, so there must be $f_{-n,0}=0$ for some $n\in\mathbb{N}(n\geq2)$, a contradiction,
and $\text {deg}_{1}f_{\pm1,0}=1$. Similarly, we can prove $\text {deg}_{2}f_{0,\pm1}=1$.
\end{proof}

\begin{lemma}
$$f_{\pm1,0}=\mu^{\pm1} (h_{1}-h_{2}\pm\nu(h_{2})), f_{0,\pm1}=e^{\pm1} (h_{2}-h_{1}\pm f(h_{1}))$$
for some $\mu, e\in\mathbb{C}^{*}$, $\nu(h_{2})\in\mathbb{C}[h_{2}], f(h_{1})\in\mathbb{C}[h_{1}]$.
\end{lemma}
\begin{proof} According to Lemma 3.2 and equation (3.2), suppose
$$f_{1,0}=\mu(h_{1}-h_{2}+\nu(h_{2})), f_{-1,0}=\mu^{-1}(h_{1}-h_{2}+\omega(h_{2}))$$
for some $\mu\in\mathbb{C}^{*}$, $\nu(h_{2}), \omega(h_{2})\in\mathbb{C}[h_{2}]$.  Then unfold (3.2)
$$\begin{aligned}
&(\sigma_{1}^{-1}f_{1,0})f_{-1,0}\\
&=(h_{1}-h_{2}+\nu(h_{2})+1)(h_{1}-h_{2}+\omega(h_{2}))\\
&=(h_{1}-h_{2})^{2}+(\nu(h_{2})+1+\omega(h_{2}))(h_{1}-h_{2})+(\nu(h_{2})+1)\omega(h_{2})\\
&=(h_{1}-h_{2})^{2}+(h_{1}-h_{2})+a_{0}(h_{2}),
\end{aligned}$$
so $\nu(h_{2})+\omega(h_{2})=0, \omega(h_{2})=-\nu(h_{2})$ and $f_{\pm1,0}=\mu^{\pm1} (h_{1}-h_{2}\pm\nu(h_{2}))$.
Similarly, we can prove $f_{0,\pm1}=e^{\pm1} (h_{2}-h_{1}\pm f(h_{1}))$.
\end{proof}

\begin{theorem}We have
$$\begin{aligned}
&f_{1,0}=\lambda_{1} (h_{1}-2h_{2}+\kappa), \\
&f_{-1,0}=\lambda_{1}^{-1} (h_{1}-\kappa),\\
&f_{0,1}=\lambda_{2}(2h_{1}-h_{2}-\kappa), \\
&f_{0,-1}=\lambda_{2}^{-1} (-h_{2}+\kappa),\\
&f_{-2,0}=\lambda_{1}^{-2}(h_{1}+h_{2}-2\kappa),\\
&f_{0,-2}=\lambda_{2}^{-2}(-h_{1}-h_{2}+2\kappa).
\end{aligned}$$
\end{theorem}

\begin{proof} By Lemma 3.3,
$$f_{-1,0}=\mu^{-1} (h_{1}-h_{2}-\nu(h_{2})), f_{0,-1}=e^{-1} (h_{2}-h_{1}- f(h_{1})).$$
Set $\tilde{\nu}(h_{2})=h_{2}+\nu(h_{2}), \tilde{f}(h_{1})=h_{1}+f(h_{1})$, then
$$f_{-1,0}=\mu^{-1} (h_{1}-\tilde{\nu}(h_{2})), f_{0,-1}=e^{-1} (h_{2}-\tilde{f}(h_{1})).$$
Consider
$$\begin{aligned}
0&=[l_{0,-1},l_{-1,0}]\cdot1=(\sigma_{2}^{-1}f_{-1,0})f_{0,-1}-(\sigma_{1}^{-1}f_{0,-1})f_{-1,0},\\
0&= (h_{1}-\tilde{\nu}(h_{2}+1))(h_{2}-\tilde{f}(h_{1}))-(h_{1}-\tilde{\nu}(h_{2}))(h_{2}-\tilde{f}(h_{1}+1)),
\end{aligned}$$
we get
\begin{equation}
\begin{split}
&(\tilde{f}(h_{1}+1)-\tilde{f}(h_{1}))h_{1}+(\tilde{\nu}(h_{2})-\tilde{\nu}(h_{2}+1))h_{2}\\
&+\tilde{f}(h_{1})\tilde{\nu}(h_{2}+1)-\tilde{f}(h_{1}+1)\tilde{\nu}(h_{2})=0.
\end{split}
\end{equation}

there must be $\text {deg}\tilde{\nu}(h_{2})=\text {deg}\tilde{f}(h_{1})=0,-1$ (i.e. $\tilde{\nu}(h_{2}), \tilde{f}(h_{1})\in\mathbb{C}$) or
$\text {deg}\tilde{\nu}(h_{2})=\text {deg}\tilde{f}(h_{1})=1$.

{\bf{Claim:}} $\tilde{\nu}(h_{2}), \tilde{f}(h_{1})\in\mathbb{C}$.

Otherwise, if $\text {deg}\tilde{\nu}(h_{2})=\text {deg}\tilde{f}(h_{1})=1$, let
$$\tilde{\nu}(h_{2})=ah_{2}+b, \tilde{f}(h_{1})=ch_{1}+d, a, c\in\mathbb{C}^{*}, b, d\in\mathbb{C},$$
then (3.3) becomes
$$(a+1)ch_{1}-a(c+1)h_{2}+ad-bc=0,$$
so $a=c=-1$ and $b=d$. Set $\kappa=b=d$ and we have
$$\begin{aligned}
&\tilde{\nu}(h_{2})=-h_{2}+\kappa, \tilde{f}(h_{1})=-h_{1}+\kappa, \kappa\in\mathbb{C},\\
&\nu(h_{2})=-2h_{2}+\kappa, f(h_{1})=-2h_{1}+\kappa, \kappa\in\mathbb{C},
\end{aligned}$$
and
$$\begin{aligned}
&f_{1,0}=\mu (h_{1}-3h_{2}+\kappa), f_{-1,0}=\mu^{-1} (h_{1}+h_{2}-\kappa),\\
&f_{0,1}=e (h_{2}-3h_{1}+\kappa),   f_{0,-1}=e^{-1} (h_{2}+h_{1}-\kappa).
\end{aligned}$$

Let's calculate
$$\begin{aligned}
f_{-1,1}&=[l_{-1,0},l_{0,1}]\cdot1\\
&=(\sigma_{1}^{-1}f_{0,1})f_{-1,0}-(\sigma_{2}f_{-1,0})f_{0,1}\\
&=e\mu^{-1}(h_{2}-3h_{1}+\kappa-3)(h_{1}+h_{2}-\kappa)\\
&-e\mu^{-1}(h_{1}+h_{2}-\kappa-1)(h_{2}-3h_{1}+\kappa)\\
&=e\mu^{-1}(-3(h_{1}+h_{2}-\kappa)+(h_{2}-3h_{1}+\kappa))\\
&=e\mu^{-1}(-6h_{1}-2h_{2}+4\kappa)\\
&=-2e\mu^{-1}(3h_{1}+h_{2}-2\kappa),
\end{aligned}$$
i.e., $f_{-1,1}=-2e\mu^{-1}(3h_{1}+h_{2}-2\kappa)$.
Check on
$$\begin{aligned}
4f_{0,1}&=[l_{-1,1},l_{1,0}]\cdot1\\
&=(\sigma_{1}^{-1}\sigma_{2}f_{1,0})f_{-1,1}-(\sigma_{1}f_{-1,1})f_{1,0}\\
&=\mu (h_{1}-3h_{2}+\kappa+4)(-2e\mu^{-1}(3h_{1}+h_{2}-2\kappa))\\
&+2e\mu^{-1}(3h_{1}+h_{2}-2\kappa-3)\mu (h_{1}-3h_{2}+\kappa)\\
&=-8e(3h_{1}+h_{2}-2\kappa)-6e(h_{1}-3h_{2}+\kappa)\\
&=e(-30h_{1}+10h_{2}+10\kappa)\\
&=10e(h_{2}-3h_{1}+\kappa),
\end{aligned}$$
i.e., $4f_{0,1}=10e(h_{2}-3h_{1}+\kappa)$ contradicts $f_{0,1}=e (h_{2}-3h_{1}+\kappa)$. The claim is true that
$\tilde{\nu}(h_{2}), \tilde{f}(h_{1})\in\mathbb{C}$. Set $\tilde{\nu}(h_{2})=\kappa, \tilde{f}(h_{1})=\iota, \kappa, \iota\in\mathbb{C}$,
then $\nu(h_{2})=-h_{2}+\kappa, f(h_{1})=-h_{1}+\iota$ and
$$\begin{aligned}
&f_{1,0}=\mu (h_{1}-2h_{2}+\kappa), f_{-1,0}=\mu^{-1} (h_{1}-\kappa),\\
&f_{0,1}=e(h_{2}-2h_{1}+\iota), f_{0,-1}=e^{-1} (h_{2}-\iota).
\end{aligned}$$
Let's calculate
$$\begin{aligned}
f_{-1,1}&=[l_{-1,0},l_{0,1}]\cdot1\\
&=(\sigma_{1}^{-1}f_{0,1})f_{-1,0}-(\sigma_{2}f_{-1,0})f_{0,1}\\
&=e\mu^{-1}(h_{2}-2h_{1}+\iota-2)(h_{1}-\kappa)-e\mu^{-1} (h_{1}-\kappa)(h_{2}-2h_{1}+\iota)\\
&=-2e\mu^{-1} (h_{1}-\kappa),\\
\end{aligned}$$
i.e., $f_{-1,1}=-2e\mu^{-1} (h_{1}-\kappa)$,
and
$$\begin{aligned}
f_{1,-1}&=[l_{1,0},l_{0,-1}]\cdot1\\
&=(\sigma_{1}f_{0,-1})f_{1,0}-(\sigma_{2}^{-1}f_{1,0})f_{0,-1}\\
&=\mu e^{-1} (h_{2}-\iota)(h_{1}-2h_{2}+\kappa)-\mu e^{-1}(h_{1}-2h_{2}+\kappa-2)(h_{2}-\iota)\\
&=2\mu e^{-1} (h_{2}-\iota),\\
\end{aligned}$$
i.e., $f_{1,-1}=2\mu e^{-1} (h_{2}-\iota)$.
By
$$\begin{aligned}
4(h_{1}-h_{2})
&=[l_{-1,1},l_{1,-1}]\cdot1\\
&=(\sigma_{1}^{-1}\sigma_{2}f_{1,-1})f_{-1,1}-(\sigma_{1}\sigma_{2}^{-1}f_{ -1,1})f_{1,-1}\\
&=-4(h_{2}-\iota-1)(h_{1}-\kappa)+4(h_{2}-\iota)(h_{1}-\kappa-1)\\
&=4(h_{1}-h_{2}-\kappa+\iota),
\end{aligned}$$
we get $-\kappa+\iota=0$, $\kappa=\iota$, and
$$\begin{aligned}
&f_{1,0}=\mu (h_{1}-2h_{2}+\kappa), f_{-1,0}=\mu^{-1} (h_{1}-\kappa),\\
&f_{0,1}=e(h_{2}-2h_{1}+\kappa), f_{0,-1}=e^{-1} (h_{2}-\kappa).
\end{aligned}$$

Now it's time to find out $f_{-2,0}$. Set $f_{-2,0}=\sum_{k=0}^{m} a_{k}(h_{2})(h_{1}-h_{2})^{k}$, consider
$$\begin{aligned}
3f_{-1,0}&=[l_{-2,0},l_{1,0}]\cdot1\\
&=(\sigma_{1}^{-2}f_{1,0})f_{-2,0}-(\sigma_{1}f_{-2,0})f_{1,0}\\
&=\mu (h_{1}-2h_{2}+\kappa+2)f_{-2,0}-(\sigma_{1}f_{-2,0})\mu (h_{1}-2h_{2}+\kappa)\\
&=\mu (h_{1}-2h_{2}+\kappa)(f_{-2,0}-(\sigma_{1}f_{-2,0}))+2\mu f_{-2,0},\\
\end{aligned}$$
which can be written out as
$$\begin{aligned}
&3\mu^{-1} (h_{1}-h_{2})+3\mu^{-1} (h_{2}-\kappa)\\
&=\mu (h_{1}-h_{2})(\sum_{k=0}^{m} a_{k}(h_{2})((h_{1}-h_{2})^{k}-(h_{1}-h_{2}-1)^{k}))\\
&+\mu (-h_{2}+\kappa)(\sum_{k=0}^{m} a_{k}(h_{2})((h_{1}-h_{2})^{k}-(h_{1}-h_{2}-1)^{k}))\\
&+2\mu \sum_{k=0}^{m} a_{k}(h_{2})(h_{1}-h_{2})^{k},
\end{aligned}$$
from which it's not difficult to see that $m=1$, and the above equation falls into
$$\begin{aligned}
&3\mu^{-1} (h_{1}-h_{2})+3\mu^{-1} (h_{2}-\kappa)\\
&=\mu (h_{1}-h_{2})a_{1}(h_{2})\\
&+\mu (-h_{2}+\kappa)a_{1}(h_{2})\\
&+2\mu (a_{0}(h_{2})+a_{1}(h_{2})(h_{1}-h_{2})),
\end{aligned}$$
which lead to
$$\begin{aligned}
&3\mu a_{1}(h_{2})=3\mu^{-1},\\
&2\mu a_{0}(h_{2})+\mu (-h_{2}+\kappa)a_{1}(h_{2})=3\mu^{-1} (h_{2}-\kappa),
\end{aligned}$$
so  $a_{1}(h_{2})=\mu^{-2}, a_{0}(h_{2})=2\mu^{-2}(h_{2}-\kappa)$ and
$$f_{-2,0}=\mu^{-2}(h_{1}+h_{2}-2\kappa).$$

Similar discussions will show that
$$f_{0,-2}=e^{-2}(-h_{1}-h_{2}+2\kappa).$$
It is done letting $\lambda_{1}=\mu$ and $\lambda_{2}=-e$.
\end{proof}

\begin{lemma}For $(r_{1},r_{2})\in\mathbb{Z}^{2}$ with $(r_{1},r_{2})\neq (-1, -1), (-2, -2)$ we have
\begin{equation}
\begin{aligned}
l_{r_{1},r_{2}}\cdot1
&=\lambda_{1}^{r_{1}}\lambda_{2}^{r_{2}}((r_{2}+1)h_{1}-(r_{1}+1)h_{2}+(r_{1}-r_{2})\kappa)\\
&=\lambda_{1}^{r_{1}}\lambda_{2}^{r_{2}}((r_{2}+1)(h_{1}-\kappa)-(r_{1}+1)(h_{2}-\kappa)).
\end{aligned}
\end{equation}
\end{lemma}

\begin{proof}
According to Theorem 3.4, this is true on the basic polynomials $f_{0,1},f_{0,-1},f_{1,0},f_{-1,0},f_{0,-2},f_{-2,0}$ corresponding to the generators,
then Definition 2 assures it.
\end{proof}

Consequently we have

\begin{theorem} Let $M$ be a $\bar{S}_{2}$-module which is a free  $U(\h_2)$-module of rank 1.
Then $M\cong \Omega(\Lambda_2,\kappa)$ for some $\Lambda_2=(\lambda_{1}, \lambda_{2})\in (\C^*)^{2}, \kappa\in \C$.
\end{theorem}

Set $$\Omega'(\Lambda_2,\kappa)=\langle h_{1}-\kappa, h_{2}-\kappa\rangle,$$
the ideal of the polynomial algebra $\C[h_{1},h_{2}]$ generated by $h_{1}-\kappa$ and $h_{2}-\kappa$.
$\Omega'(\Lambda_2,\kappa)$ is a submodule considering (3.4). For the simplicity, we provide

\begin{theorem}$\Omega'(\Lambda_2,\kappa)$ is the unique irreducible submodule of $\Omega(\Lambda_2,\kappa)$, and
$$\Omega(\Lambda_2,\kappa)/\Omega'(\Lambda_2,\kappa)\cong\mathbb{C}.$$
\end{theorem}

\begin{proof} The isomorphism is immediate. For the irreducibility, take any nonzero $F(h_{1},h_{2})\in\Omega(\Lambda_2,\kappa)$, let $N=U(\bar{S}_{2})(F(h_{1},h_{2}))$, it's sufficient to prove $h_{1}-\kappa,h_{2}-\kappa\in N$. Observing
$$l_{i,-1}\cdot F(h_{1},h_{2})=-(i+1)\lambda_{1}^{i}\lambda_{2}^{-1}F(h_{1}-i,h_{2}+1)(h_{2}-\kappa), i\in\mathbb{Z},$$
 we see that there is a nonzero $G(h_{2})\in N\cap\mathbb{C}[h_{2}]$. Considering
$$l_{-1,j}\cdot G(h_{2})=(j+1)\lambda_{1}^{-1}\lambda_{2}^{j}G(h_{2}-j)(h_{1}-\kappa), j\in\mathbb{Z},$$
we see that $h_{1}-\kappa\in N$.
Besides,
$$h_{2}-\kappa=\frac{1}{2}\lambda_{1}^{-1}\lambda_{2}l_{1,-1}\cdot(h_{1}-\kappa)+(h_{1}-\kappa)(h_{2}-\kappa)\in N.$$
We also see that $\Omega'(\Lambda_2,\kappa)$ is the unique nontrivial submodule, the uniqueness is clear by the proof of the irreducibility..
\end{proof}

The simple module $\Omega'(\Lambda_2,\kappa)$ is completely determined by the triple $(\Lambda_2,\kappa)=(\lambda_{1},\lambda_{2},\kappa)$.
\begin{theorem}Let $\Lambda_2=(\lambda_{1}, \lambda_{2}), \Lambda_2'=(\lambda_{1}',\lambda_{2}')\in (\C^*)^{2}, \kappa, \kappa'\in \C$. Then the simple $\bar{S}_{2}$-modules
$\Omega'(\Lambda_2,\kappa)\cong\Omega'(\Lambda_2',\kappa')$ if and only if
$\lambda_{1}=\lambda_{1}',\lambda_{2}=\lambda_{2}',\kappa=\kappa'$.
\end{theorem}

\begin{proof} The sufficiency is evident. For the necessity, let $$\psi:\Omega'(\Lambda_2,\kappa))\rightarrow\Omega'(\Lambda_2',\kappa')$$
be an isomorphism between the two $\bar{S}_{2}$-modules. Let's consider
$$(l_{-1,0}-\lambda_{1}^{-1}(h_{1}-\kappa))(h_{2}-\kappa)=0,$$
letting $\psi$ act on it and we get
$$\begin{aligned}
&(l_{-1,0}-\lambda_{1}^{-1}(h_{1}-\kappa))\psi(h_{2}-\kappa)=0,\\
&\lambda_{1}'^{-1}(\sigma_{1}^{-1}(\psi(h_{2}-\kappa)))(h_{1}-\kappa')=\lambda_{1}^{-1}(h_{1}-\kappa)\psi(h_{2}-\kappa),
\end{aligned}$$
so $\lambda_{1}=\lambda_{1}'$. Rewrite the previous equality
$$(\sigma_{1}^{-1}(\psi(h_{2}-\kappa))-\psi(h_{2}-\kappa))(h_{1}-\kappa')=(\kappa'-\kappa)\psi(h_{2}-\kappa),$$
we see that $\kappa'-\kappa\in\mathbb{Z_{+}}$, and  $\kappa-\kappa'\in\mathbb{Z_{+}}$ by symmetry, so $\kappa'=\kappa$, and $\psi(h_{2}-\kappa)\in\mathbb{C}[h_{2}]$. By similar discussions, we can show that $\psi(h_{1}-\kappa)\in\mathbb{C}[h_{1}]$.
\end{proof}

\subsection{$\bar{S}_{n}$- module structures on $M=U(\h_n)(n\geq3)$}

Now it's time to work from $2$ to $n (n\geq3)$. Let's consider the $\bar{S}_{n}$- module structures on $M=U(\h_n)\cdot1\cong\mathbb{C}[\partial_1, \partial_2, \cdot\cdot\cdot, \partial_n]$ that is a free $U(\h_n)$-module of rank 1, on which $\bar{S}_{n}'=t^{-1}\tilde{S}_{n}'$ acts nontrivially,
where, $\mathfrak{h}_{n}=\oplus_{i=1}^n\mathbb{C} \partial_i$.

Note that $\bar{S}_{n}'\cap\mathfrak{h}_{n}=\bigoplus_{i=1}^{n-1}\mathbb{C} (\partial_i-\partial_{i+1})$  and $\bar{S}_{n}'$ is simple, so a free $U(\h_{n})$-module of rank 1 is naturally a module on which $\bar{S}_{n}'$ acts nontrivially.

\begin{lemma} For any $f(\partial_1,\partial_2,...,\partial_n)\in M$, $r=(r_1,r_2,\cdot\cdot\cdot,r_n)\in \mathbb{Z}^{n}$,
$i, j\in\{1,\cdot\cdot\cdot,n\}, i\neq j$, we have
$$\begin{aligned}
&((r_{j}+1)t^{r}\partial_i-(r_{i}+1)t^{r}\partial_j)\cdot f(\partial_1,\partial_2,...,\partial_n)1\\
&=f(\partial_1-r_1,\partial_2-r_2,...,\partial_n-r_n)(((r_{j}+1)t^{r}\partial_i-(r_{i}+1)t^{r}\partial_j)\cdot 1),\\
\end{aligned}$$
and $((r_{j}+1)t^{r}\partial_i-(r_{i}+1)t^{r}\partial_j)\cdot 1\ne0$.
\end{lemma}

\begin{proof} The proof is similar to that of Lemma 3.1.
\end{proof}

Note that $\bar{S}_{n}$ admits
$\left(
\begin{array}{c}
n\\
2\\
\end{array}
\right)$ subalgebras:
$$\bar{S}_{ij}=Span_{\mathbb{C}}\{\partial_{i}, \partial_{j}, l_{r_{i},r_{j}}\mid (r_{i},r_{j})\in \mathbb{Z}^{2}, (r_{i},r_{j})\neq(-1,-1),(-2,-2)\},$$
and
$$\bar{S}_{ij}'=Span_{\mathbb{C}}\{l_{r_{i},r_{j}}\mid (r_{i},r_{j})\in \mathbb{Z}^{2}, (r_{i},r_{j})\neq(-1,-1),(-2,-2)\},$$
where, $l_{r_{i},r_{j}}=(r_{j}+1)t_{i}^{r_{i}}t_{j}^{r_{j}}\partial_i-(r_{i}+1)t_{i}^{r_{i}}t_{j}^{r_{j}}\partial_j ,i, j\in\{1,\cdot\cdot\cdot,n\}, i\neq j$.

It's immediate that $\bar{S}_{ij}=\bar{S}_{ji}$ and $\bar{S}_{ij}\cong \bar{S}_{2}$.

\begin{lemma}
$\bar{S}_{ij}-$actions on $M$ are compatible.
\end{lemma}

\begin{proof}
In the following, $i, j, l\in\{1,\cdot\cdot\cdot,n\}$, which are pairwise different.
In  $\bar{S}_{ij}$,
$$l_{r_{i},r_{j}}\cdot1=\lambda_{ij}^{r_{i}}\lambda_{ji}^{r_{j}}((r_{j}+1)\partial_{i}-(r_{i}+1)\partial_{j}+(r_{i}-r_{j})\kappa_{ij}),$$
and
$$t_{i}^{-1}\partial_{i}\cdot1=l_{-1,0}\cdot1=\lambda_{ij}^{-1}(\partial_{i}-\kappa_{ij}), \kappa_{ij}\in R_{ij};$$
while in $\bar{S}_{il}$,
$$t_{i}^{-1}\partial_{i}\cdot1=\lambda_{il}^{-1}(\partial_{i}-\kappa_{il}), \kappa_{il}\in R_{il}.$$
We get $\lambda_{ij}^{-1}(\partial_{i}-\kappa_{ij})=\lambda_{il}^{-1}(\partial_{i}-\kappa_{il})$, so $\lambda_{ij}=\lambda_{il}$ and $\kappa_{ij}=\kappa_{il}\in\mathbb{C}$, the compatibility is established.
\end{proof}

\begin{lemma}
$\left\{\bar{S}_{ij}\mid i, j\in\{1,\cdot\cdot\cdot,n\} ,i\neq j \right\}$ generate $\bar{S}_{n}$.
\end{lemma}

\begin{proof}
Define the length
$$\ell((r_{j}+1)t^{r}\partial_i-(r_{i}+1)t^{r}\partial_j)=\ell(t^{r})=\lvert \{r_{i}\neq0,i=1,\cdot\cdot\cdot,n\}\rvert.$$
Make induction on $\ell(t^{r})$. For pairwise different $i, j, l\in\{1,\cdot\cdot\cdot,n\}$, observe
$$\begin{aligned}
&[(r_{j}+1)t_{i}^{r_{i}}t_{j}^{r_{j}}\partial_i-(r_{i}+1)t_{i}^{r_{i}}t_{j}^{r_{j}}\partial_j, (r_{l}+1)t_{l}^{r_{l}}\partial_{i}-t_{l}^{r_{l}}\partial_{l}]\\
&=-(r_{l}+1)r_{i}((r_{j}+1)t_{i}^{r_{i}}t_{j}^{r_{j}}t_{l}^{r_{l}}\partial_i-(r_{i}+1)t_{i}^{r_{i}}t_{j}^{r_{j}}t_{l}^{r_{l}}\partial_j),
\end{aligned}$$
and
$$\begin{aligned}
&[(r_{j}+1)t_{i}^{r_{i}}t_{j}^{r_{j}}t_{l}^{-2}\partial_i-(r_{i}+1)t_{i}^{r_{i}}t_{j}^{r_{j}}t_{l}^{-2}\partial_j, 2t_{l}\partial_{i}-t_{l}\partial_{l}]\\
&=-2(r_{i}+1)((r_{j}+1)t_{i}^{r_{i}}t_{j}^{r_{j}}t_{l}^{-1}\partial_i-(r_{i}+1)t_{i}^{r_{i}}t_{j}^{r_{j}}t_{l}^{-1}\partial_j),
\end{aligned}$$
so $(r_{j}+1)t^{r}\partial_i-(r_{i}+1)t^{r}\partial_j$ can be generated when $\ell(t^{r})=3$. Suppose
$(r_{j}+1)t_{1}^{r_{1}}\cdot\cdot\cdot \overset{\wedge}{t_{s}^{r_{s}}}\cdot\cdot\cdot t_{n}^{r_{n}}\partial_i-(r_{i}+1)t_{1}^{r_{1}}\cdot\cdot\cdot \overset{\wedge}{t_{s}^{r_{s}}}\cdot\cdot\cdot  t_{n}^{r_{n}}\partial_j$ with
$\ell(t_{1}^{r_{1}}\cdot\cdot\cdot \overset{\wedge}{t_{s}^{r_{s}}}\cdot\cdot\cdot t_{n}^{r_{n}})=n-1(n\geq4)$ can be generated (note that $i, j, s\in\{1,\cdot\cdot\cdot,n\}$ are pairwise different), then
$$\begin{aligned}
&[(r_{j}+1)t_{1}^{r_{1}}\cdot\cdot\cdot \overset{\wedge}{t_{s}^{r_{s}}}\cdot\cdot\cdot t_{n}^{r_{n}}\partial_i-(r_{i}+1)t_{1}^{r_{1}}\cdot\cdot\cdot \overset{\wedge}{t_{s}^{r_{s}}}\cdot\cdot\cdot  t_{n}^{r_{n}}\partial_j, (r_{s}+1)t_{s}^{r_{s}}\partial_{i}-t_{s}^{r_{s}}\partial_{s}]\\
&=-(r_{s}+1)r_{i}((r_{j}+1)t^{r}\partial_i-(r_{i}+1)t^{r}\partial_j),
\end{aligned}$$
and
$$\begin{aligned}
&[(r_{j}+1)t_{1}^{r_{1}}\cdot\cdot\cdot t_{s}^{-2}\cdot\cdot\cdot t_{n}^{r_{n}}\partial_i-(r_{i}+1)t_{1}^{r_{1}}\cdot\cdot\cdot {t_{s}^{-2}}\cdot\cdot\cdot  t_{n}^{r_{n}}\partial_j, 2t_{s}\partial_{i}-t_{s}\partial_{s}]\\
&=-2(r_{i}+1)((r_{j}+1)t_{1}^{r_{1}}\cdot\cdot\cdot t_{s}^{-1}\cdot\cdot\cdot t_{n}^{r_{n}}\partial_i-(r_{i}+1)t_{1}^{r_{1}}\cdot\cdot\cdot {t_{s}^{-1}}\cdot\cdot\cdot  t_{n}^{r_{n}}\partial_j),
\end{aligned}$$
so $(r_{j}+1)t^{r}\partial_i-(r_{i}+1)t^{r}\partial_j$ of length $n$ can be generated.

Besides, it's inviting to see that $l_{-2,-2}=-t_{i}^{-2}t_{j}^{-2}\partial_i+t_{i}^{-2}t_{j}^{-2}\partial_j$ can be produced (which can never be achieved in a single $\bar{S}_{ij}$). Take $t_{i}^{-1}t_{j}^{-1}\partial_i$  and $t_{i}^{-1}t_{j}^{-1}\partial_j$, they make it:
$$[t_{i}^{-1}t_{j}^{-1}\partial_i,t_{i}^{-1}t_{j}^{-1}\partial_j]=t_{i}^{-2}t_{j}^{-2}\partial_i-t_{i}^{-2}t_{j}^{-2}\partial_j.$$
\end{proof}

The $\bar{S}_{n}$-module structures on $M=U(\h_n)\cdot1$ are:
\begin{theorem} Let $M$ be a $\bar{S}_{n}$-module which is a free $U(\h_n)$-module of rank 1.
Then $M\cong\Omega(\Lambda_n,\kappa)$ for some $\Lambda_n\in (\C^*)^n, \kappa\in \C$.
\end{theorem}
\begin{proof} This follows from Theorem 3.6, Lemma 3.9, Lemma 3.10 and Lemma 3.11.
\end{proof}

Set $$\Omega'(\Lambda_n,\kappa)=\langle\partial_1-\kappa, \cdot\cdot\cdot, \partial_n-\kappa\rangle,$$
the ideal of the polynomial algebra $\mathbb{C}[\partial_1, \partial_2, \cdot\cdot\cdot, \partial_n]$ generated by $\{\partial_i-\kappa\mid i=1\cdot\cdot\cdot, n\}$. $\Omega'(\Lambda_n,\kappa)$ is a submodule of $\Omega(\Lambda_n,\kappa)$. By similar discussions, we can prove that
\begin{theorem}Let $\Lambda_n\in (\C^*)^n, \kappa\in \C$. Then
$\Omega(\Lambda_n,\kappa)$ has a unique irreducible submodule $\Omega'(\Lambda_n,\kappa)$,
and the quotient is $\mathbb{C}$.
\end{theorem}

\begin{theorem}Let $\Lambda_n, \Lambda_n'\in (\C^*)^n, \kappa, \kappa'\in\C$. Then
the simple $\bar{S}_{n}$-modules
$\Omega'(\Lambda_n,\kappa)\cong \Omega'(\Lambda'_n,\kappa')$ if and only if $\Lambda_n=\Lambda'_n,\kappa=\kappa'$.
\end{theorem}

\noindent {\bf Acknowledgement.}

\vskip 5pt

Sincere thanks should go to Professor Kaiming Zhao for formulating the problem and stimulating discussions,
for his keen guidance and good advices on the paper. Special thanks should also go to Professor Xiandong Wang for his patient discussions.
The author is partially supported by NSFC (Grant No.11472144)
and NSFC (Grant No.11501315).

\vspace{10mm}

\end{document}